\documentclass[3p,times]{elsarticle}
\usepackage{amsfonts}
\usepackage{amssymb}
\usepackage{amsmath}
\usepackage{fancyvrb}
\usepackage{graphics}
\usepackage{graphicx}
\usepackage{epsfig}
\usepackage{epstopdf}
\usepackage{lipsum}
\usepackage{leftidx}
\usepackage{multirow}
\usepackage{charter}
\usepackage[charter]{mathdesign}
\usepackage{color}
\usepackage{amsmath}
\usepackage{ntheorem}
\usepackage{float}
\usepackage{subfig}
\newtheorem*{proof}{Proof}

\begin{document}
\newtheorem{thm}{Theorem}
\newtheorem{lem}{Lemma}
\begin{frontmatter}
\title{High Order Compact Finite Difference Methods for Non-Fickian Flows in Porous Media}
   \author[label1]{Xuan Zhao}
   \author[label2]{Ziyan Li}
   \author[label1]{Xiaoli Li\corref{cor1}}
   \cortext[cor1]{Corresponding author.}
   \ead{xiaolimath@sdu.edu.cn.}
   \address[label1]{School of Mathematics, Shandong University, Jinan, Shandong, 250100, China.}
   \address[label2]{Department of Mathematics, City University of Hong Kong, Hong Kong SAR, China.}
\begin{abstract}
    In this work,  fourth-order compact block-centered finite difference (CBCFD) schemes combined with the Crank-Nicolson discretization are constructed and analyzed for solving parabolic integro-differential type non-Fickian flows in one-dimensional and two-dimensional cases. Stability analyses of the constructed schemes are derived rigorously. We also obtain the optimal second-order convergence in temporal increment and the fourth-order convergence in spatial direction for both velocity and pressure. To verify the validity of the CBCFD schemes, we present some experiments to show that the numerical results are in agreement with our theoretical analysis.   

\end{abstract}


\end{frontmatter}
 \thispagestyle{empty}
\section{Introduction}
The non-Fickian flow is widely used to describe the transport of contaminants in porous media, which is complicated by the history effect which characterizes various mixing length growth of the flow. Actually the evolution of a reactive chemical within a velocity field usually represents by using the classical Fickian dispersion theory. For instance, the evolution in such a velocity field, when modeled with Fickian-type constitutive laws, leads to a dispersion tensor dependent upon the timescales of observation. Hence, to avoid this difficulty, non-Fickian models have been recently proposed, in which the dispersion term arising from integration with respect to time makes the flow non-Fickian, since it is not a pure diffusion term \cite{ewing2002sharp}.  The parabolic integro-differential model is a type of non-Fickian model that can represent many physical processes such as non-local groundwater transport \cite{Dagan1994The}, microsensor thermistor problems \cite{1999ALin} and so on.

Due to the exists of various mixing length growth, non-Fickian flows are complicated and also difficult to be solved \cite{1988Non,2005L}. Therefore, it is meaningful for us to pay close attention to the numerical methods for solving such models. There is sizeable literature on the numerical approximations of the problem.  In 2000, Ewing \cite{2000Finite} constructed various finite volume element schemes for the parabolic integro-differential problem and demonstrated the error estimates in $H^{1}$-norm and $L^{2}$-norm.  Arikoglu and Ozkol \cite{arikoglu2005solution} considered the differential transform method for solving the integro-differential equations and introduced new theorems to show the fast convergence of the method. Rui and Guo \cite{2013Split} established a least-squares finite element method and demonstrated the optimal convergence order.  Besides, some other methods were proposed such as mixed finite element method \cite{2002LJ}, two-grid method \cite{2014A} and so on.

Due to its simplicity and high efficiency, the finite difference method is widely used in practical engineering projects. Furthermore, the block-centered finite difference (BCFD) method based on the staggered grids can maintain local mass conservation, which is of great significance for engineering calculation. The BCFD method which can be thought of as a type of mixed element method \cite{raviart_1977} has been applied to solve many types of partial differential equations. Wheeler and Weiser \cite{1988On} derived the discrete $L^{2}$-norm error estimate of the BCFD method for the linear elliptic equations and proved the second-order convergence on nonuniform grids. In 2002, the BCFD method with local grid refinement was applied to solve the groundwater problem \cite{2002Development}. The authors demonstrated the validity for both the homogenous and heterogeneous systems. Rui and Pan \cite{2012Pan,2015A} introduced the BCFD scheme to solve the nonlinear Darcy-Forchheimer model with constant coefficient and variable Forchheimer number and proved the second-order convergence for both pressure and velocity. In 2015, the combination of the two-grid and BCFD method was constructed for solving the Darcy-Forchheimer model and the optimal convergence order in discrete $L^{2}$-norm was obtained \cite{2015Rui}. Some other works on different models can be found in \cite{2015Zhai,2015Zh,2016liuz,liu2018block,li2018block}. In 2018, Li and Rui proposed the BCFD method for the non-Fickian flow model \cite{li2018blocknoF}. They introduced two schemes with the difference in temporal discretization and analyzed the convergences rigorously. This was the first attempt to solve the non-Fickian flow model by using the BCFD method. However, the classical BCFD method can only obtain second-order  convergence in the spatial direction, which may be not enough in some high-accuracy cases. Therefore, the construction of the high order numerical scheme is of great significance. In 2021, Xie et al. \cite{2021High} proposed fourth-order CBCFD schemes for      
general elliptic and parabolic problems. The authors constructed the high order operators, analyzed one-dimensional and two-dimensional problems, and obtained the fourth order convergence for both pressure and velocity in the spatial direction. Based on the work, the high order scheme was also applied to solve the nonlinear contaminant transport model \cite{SHI2022212}. 

As far as we know, there is no research using the CBCFD scheme to solve the non-Fickian flows in the porous media. Therefore in this work, we combine the compact block-centered finite difference method in spatial discretization and the Crank-Nicolson scheme in temporal discretization to solve non-Fickian flows in porous media, which can lead to a higher order convergence in spatial direction compared with the classical block-centered finite difference method. We establish the stability analyses and error estimates of the constructed schemes rigorously both in one-dimensional and two-dimensional cases. Compared with theoretical analysis in one-dimensional problem, the main difficulty for two-dimensional case is that we should establish new boundary preserving properties since the higher-order operators acting on two directions. And also compared with the recent works in \cite{2021High,SHI2022212}, the main contribution in this paper lies in the careful and special treatment of the nonlocal term to obtain the optimal convergence. Finally, we present some numerical experiments to verify the theoretical analysis. 

This paper is organized as follows. In Sect. 2, we introduce the one-dimensional problem and derived the stability and convergence. In Sect. 3, we consider the two-dimensional problem based on the one-dimensional case. In Sect. 4, some numerical experiments are carried out. 

Throughout this whole paper, we use $\epsilon$ and $C$, with or without subscript, to denote a small positive constant and a positive constant respectively, which may have different values at different appearances.
\biboptions{numbers,sort&compress}

\section{One-dimensional parabolic problem}
In this section, the following one-dimensional non-Fickian flow in porous media is considered: 
\begin{equation}\label{e1}
  \left\{\begin{array}{ll} p_{t}+u_{x}=f(x,t), \quad (x, t) \in \Omega\times J, \\
  u=-a\left(x\right)p_{x}\left(x,t\right)-\int_{0}^{t}b\left(x,s\right)p_{x}\left(x,s\right)ds, \quad (x, t) \in \Omega\times J,\\
  -a\left(x\right)p_{x}\left(x,t\right)=0, \quad x \in \partial \Omega, t \in  J,\\
  p\left(x,0\right)=p_{0}(x),  \quad x \in \Omega. 
\end{array}\right.
\end{equation}
Here $\Omega$ is a one-dimensional domain defined as $\left(0,L\right)$ and $J$ is $\left(0,T\right]$.  We suppose $a(x),$ $a^{'}(x),$ $b(x,t)$ and $b^{'}(x,t)$ are bounded smooth functions and satisfy
\begin{equation}\label{e2}
  \begin{split}
  a_{*}\leq a\left(x\right)\leq a^{*},\quad b_{*}\leq b\left(x,t\right)\leq b^{*},
  \end{split}
\end{equation}
where $a_{*},a^{*},b_{*}$ and $ b^{*}$ are positive constants.

 Set $\widetilde{u}=-a\left(x\right)p_{x}\left(x,t\right)$.  The Eq.$(\ref{e1})$ can be recast into the following formulation:
 \begin{equation}\label{e3}
  \left\{\begin{array}{ll} p_{t}+u_{x}=f(x,t), \quad (x, t) \in \Omega\times J, \\
  u=\widetilde{u}+\int_{0}^{t}\frac{b}{a}\widetilde{u}\left(x,s\right)ds, \quad (x, t) \in \Omega\times J,\\
  \widetilde{u}\left(x,t\right)=0, \quad x \in \partial \Omega, t \in  J,\\
  p\left(x,0\right)=p_{0}(x),  \quad x \in \Omega. 
\end{array}\right.
\end{equation}
\subsection{The CBCFD scheme of one dimensional problem}
First, we give some notations which will be helpful for the analysis of the CBCFD scheme.

Define $\Omega_{i}=[x_{i-1/2},x_{i+1/2}]$, $i=1,2,\cdots, M $ as the uniform partitions of $\Omega$, where $x_{1/2}=0,$ and $x_{N+1/2}=L.$ Set $x_{i}=(x_{i-1/2}+x_{i+1/2})/2,$ $h=L/M,$ and $x_{i+1/2}=ih.$ Let $f$ be any function, we give the definitions of the following difference operators
\begin{equation*}
  \begin{split}
    \delta_{x}f_{i}&=\frac{f_{i+1/2}-f_{i-1/2}}{h},\\
    \delta_{x}^{2}f_{i}&=\frac{f_{i+1}-2f_{i}+f_{i-1}}{h^{2}},\\
    \delta_{x}f_{i+1/2}&=\frac{f_{i+1}-f_{i}}{h},\\
    \delta_{x}^{2}f_{i+1/2}&=\frac{f_{i+3/2}-2f_{i+1/2}+f_{i-1/2}}{h^{2}}.
  \end{split}
\end{equation*}
Supposing $f'(x)=g(x),$ then we have 
\begin{equation*}
  \begin{split}
    \delta_{x}f_{i}&=\frac{g_{i+1}+22g_{i}+g_{i-1}}{24}+O(h^{4}),\\
    \delta_{x}f_{i+1/2}&=\frac{g_{i+3/2}+22g_{i+1/2}+g_{i-1/2}}{24}+O(h^{4}).
  \end{split}
\end{equation*}
Now we define the interpolation operators as follows \cite{2021High},
\begin{equation*}
  \begin{split}
    \psi_{x}g_{i+1/2}&=\frac{g_{i+3/2}+22g_{i+1/2}+g_{i-1/2}}{24}\\
    &=\left(I+\frac{h^{2}}{24}\delta_{x}^{2}\right)g_{i+1/2},\qquad i=1,2,\cdots, M-1,
  \end{split}
\end{equation*}
\begin{equation*}
  \widetilde{\psi}_{x} g_{i}= \left\{\begin{aligned}
  &\left(I+\frac{h^{2}}{24} \delta_{x}^{2}\right) g_{i}, \quad i=2,3, \cdots, M-1, \\
   &  \frac{g_{1 / 2}+4 g_{1}+g_{3 / 2}}{6},\qquad\quad i=1, \\
   &   \frac{g_{J-1 / 2}+4 g_{J}+g_{J+1 / 2}}{6},  \quad i=M,
    \end{aligned}
    \right.
  \end{equation*}
  \begin{equation*}
    \widehat{\psi}_{x} g_{i}= \left\{\begin{aligned}
    &\left(I+\frac{h^{2}}{24} \delta_{x}^{2}\right) g_{i}, \quad i=2,3, \cdots, M-1, \\
     &  \frac{26g_{1}-5g_{2}+4g_{3}-g_{4}}{6},\qquad\quad i=1, \\
     &   \frac{26g_{J}-5 g_{J-1}+4g_{J-2}-3g_{J-3}}{24},  \quad i=M.
      \end{aligned}
      \right.
    \end{equation*}
Supposing $u_{t}\in C^{2}(0,1)$, then by straightforward calculation, we have 
\begin{equation*}
  \int_{0}^{t_{n+1 / 2}} u(t) d t=\Delta t \sum_{l=0}^{n-1} u(t_{l+1 / 2})+\frac{\Delta t}{2} u(t_{n+1 / 2})+O(\Delta t^{2}).
  \end{equation*}
Let $P^{n},$ $U^{n}$ and $\widetilde{U}^{n}$ be the approximations to $p^{n},$ $u^{n}$ and $\widetilde{u}^{n},$ respectively. We now construct the CBCFD scheme for one-dimensional non-Fickian flow model, 
\begin{equation}\label{e4}
  \begin{split}
    \frac{\widehat{\psi}_{x}P_{i}^{n+1}-\widehat{\psi}_{x}P_{i}^{n}}{\Delta t}+\delta_{x}U_{i}^{n+\frac{1}{2}}=\widetilde{\psi}_{x}f_{i}^{n+\frac{1}{2}},\quad i=1,2,\cdots,M,
  \end{split}
\end{equation}
\begin{equation}\label{e5}
  \begin{split}
    \delta_{x} P_{i+\frac{1}{2}}^{n+1}=-\psi_{x}\left(\frac{\widetilde{U}}{a}\right)_{i+\frac{1}{2}}^{n+1},\quad i=1,2,\cdots,M-1,
  \end{split}
\end{equation}
\begin{equation}\label{e6}
  \begin{split}
    U_{i+\frac{1}{2}}^{n+1}=\widetilde{U}_{i+\frac{1}{2}}^{n+1}+\Delta t\sum_{l=0}^{n}\left(\frac{b}{a}\widetilde{U}\right)_{i+\frac{1}{2}}^{l+\frac{1}{2}},\quad i=1,2,\cdots,M-1.
  \end{split}
\end{equation}
Set the boundary conditions and initial approximation as follows:
$$U_{\frac{1}{2}}^{n}=\widetilde{U}_{\frac{1}{2}}^{n}=0,\quad U_{M+\frac{1}{2}}^{n}=\widetilde{U}_{M+\frac{1}{2}}^{n}=0,\quad P_{i}^{0}=p_{0}(x_{i}).$$
To derive stability and convergence, we first give some Lemmas. 
\begin{lem}\label{lem1}\cite{wheeler1988}
  Suppose $q_{i}$ and $w_{i+1/2}$  are any values such that
  $w_{1/2}=w_{M+1/2}=0$, then
    \begin{equation*}
      \begin{split}
      &\left(q, \delta_{x} w\right)=-\left(\delta_{x} q, w\right).
      \end{split}
      \end{equation*}
    \end{lem}
By following exactly the same procedure in \cite{2021High}, we can obtain the following three Lemmas.
\begin{lem}\label{addlem1}\cite{2021High}
  Suppose $V\in U_{h}, a(x)\in C^{1}(\Omega)$ and $C_{a}:=\frac{\left\lVert a'(x)\right\rVert_{\infty} }{24a_{*}}$, we have
    \begin{equation*}
      \begin{split}
      (aV,\psi_{x}V)\geq\left(\frac{5}{6}-hC_{a}\right)(aV,V).
      \end{split}
      \end{equation*}
    \end{lem}
    \begin{lem}\label{addlem2}\cite{2021High}
     Set $D_{a}:=\left\lVert (1/a)'\right\rVert_{\infty} $ and $C_{a}:=\frac{\left\lVert a'(x)\right\rVert_{\infty} }{24a_{*}}.$ 
     Let $\Delta t$ be the time step and $\epsilon$ be a small constant, we obtain the following estimate
        \begin{equation*}
          \begin{split}
          \sum_{n=1}^{N-1}\left(\widetilde{U}^{n+\frac{1}{2}},\psi_{x}\frac{\partial_{t}\widetilde{U}^{n+1}}{a}\right)\geq &\frac{1}{2\Delta t}\left[\left(\frac{5}{6a^{*}}-\left(\frac{C_{a}}{a^{*}}+\frac{D_{a}}{24}\right)h\right)\left\lVert \widetilde{U}^{N}\right\rVert^{2}-\left(\frac{1}{a_{*}}+\frac{D_{a}h}{24}\right)\left\lVert \widetilde{U}^{0}\right\rVert^{2}\right.\\   
          &\left.-\frac{8\left(a^{*}\right)^{2}\epsilon}{33}\sum_{n=0}^{N-1}\Delta t\left\lVert \partial_{t}P^{n+1}\right\rVert^{2}-\frac{D_{a}^{2}}{24\epsilon}\sum_{n=0}^{N}\Delta t\left\lVert \widetilde{U}^{n}\right\rVert^{2}\right].
          \end{split}
          \end{equation*}
    \end{lem}
    \begin{lem}\label{addlem3}\cite{2021High}
      Suppose $P$ and $U$ are the numerical solutions of the CBCFD scheme, we have 
         \begin{equation*}
           \begin{split}
            \frac{3}{4}\left\lVert P\right\rVert^{2}&\leq(\widehat{\psi}_{x}P,P)\leq \frac{4}{3}\left\lVert P\right\rVert^{2},\quad (\widehat{\psi}_{x}P,\widehat{\psi}_{x}P)\leq \frac{5}{3}\left\lVert P\right\rVert^{2},\\
            \frac{5}{6}\left\lVert U\right\rVert^{2}&\leq (\psi_{x}U,U)\leq \left\lVert U\right\rVert^{2}, \quad \frac{11}{16}\left\lVert U\right\rVert^{2}\leq (\psi_{x}U,\psi_{x}U)\leq \left\lVert U\right\rVert^{2},  
           \end{split}
           \end{equation*}
           \begin{equation*}
            \begin{split}
             \hspace{-4.4cm}\lVert \psi_{x}\frac{\partial_{t}U^{n+1}}{a}\rVert^{2}\leq \frac{4}{h^{2}}\lVert \partial_{t}P^{n+1}\rVert^{2}.  
            \end{split}
            \end{equation*}
     \end{lem}
\subsection{The analysis of stability }
Multiplying Eq.$(\ref{e4})$ by $h_{i}P_{i}^{n+\frac{1}{2}},$ summing on $i$ and using Lemma $\ref{lem1}$, we obtain
\begin{equation}\label{e7}
  \left(\partial_{t} \widehat{\psi}_{x} P^{n+1}, P^{n+\frac{1}{2}}\right)+\left(U^{n+\frac{1}{2}}, \psi_{x}\left(\frac{\widetilde{U}}{a}\right)^{n+\frac{1}{2}}\right)=\left(\widetilde{\psi}_{x} f^{n+\frac{1}{2}}, P^{n+\frac{1}{2}}\right).
  \end{equation}
By simple calculation, we can get the following equation easily,
\begin{equation}\label{e8}
  \begin{array}{r}
  \left(\partial_{t} \widehat{\psi}_{x} P^{n+1}, P^{n+\frac{1}{2}}\right)=\frac{1}{\Delta t}\left(\left(\widehat{\psi}_{x} P^{n+1}, P^{n+1}\right)-\left(\widehat{\psi}_{x} P^{n}, P^{n}\right)\right)-\left(\widehat{\psi}_{x} P^{n+\frac{1}{2}}, \partial_{t} P^{n+1}\right).
  \end{array}
  \end{equation}
Using the Eq.$(\ref{e6})$, $U^{n+\frac{1}{2}}$ can be expressed as
\begin{equation*}
  \begin{split}
  U^{n+\frac{1}{2}}&=\frac{U^{n+1}+U^{n}}{2}\\
  &=\widetilde{U}^{n+\frac{1}{2}}+\Delta t \sum_{l=0}^{n-1}\left(\frac{b}{a}\widetilde{U}\right)^{l+\frac{1}{2}}+\frac{1}{2}\Delta t\left(\frac{b}{a}\widetilde{U}\right)^{n+\frac{1}{2}}. 
  \end{split}
  \end{equation*}
Therefore,
\begin{equation}\label{e9}
  \begin{split}
  \left(U^{n+\frac{1}{2}}, \psi_{x}\left(\frac{\widetilde{U}}{a}\right)^{n+\frac{1}{2}}\right)=&\left(\widetilde{U}^{n+\frac{1}{2}},\psi_{x}\left(\frac{\widetilde{U}}{a}\right)^{n+\frac{1}{2}}\right)+\left(\Delta t \sum_{l=0}^{n-1}\left(\frac{b}{a}\widetilde{U}\right)^{l+\frac{1}{2}},\psi_{x}\left(\frac{\widetilde{U}}{a}\right)^{n+\frac{1}{2}}\right)\\
  &+\left(\frac{1}{2}\Delta t\left(\frac{b}{a}\widetilde{U}\right)^{n+\frac{1}{2}},\psi_{x}\left(\frac{\widetilde{U}}{a}\right)^{n+\frac{1}{2}}\right)\\
  &=T_{1}+T_{2}+T_{3}.
  \end{split}
\end{equation}
By the Lemma $\ref{addlem1}$, the estimate of $T_{1}$ can be easily obtained.
\begin{equation}\label{e10}
  \begin{split}
   T_{1}\geq \frac{5-6C_{a}h}{6a^{*}}\left\lVert \widetilde{U}^{n+\frac{1}{2}}\right\rVert^{2}. 
  \end{split}
\end{equation}
Use the Cauchy-Schwarz inequality and the Lemma $\ref{addlem3}$, we have that 
\begin{equation}\label{e11}
  \begin{split}
   T_{2}&\geq -C\Delta t\sum_{l=0}^{n-1}\left\lVert \left(\frac{b}{a}\widetilde{U}\right)^{l+\frac{1}{2}}\right\rVert^{2}- \epsilon\Delta t\left\lVert \psi_{x}\left(\frac{\widetilde{U}}{a}\right)^{n+\frac{1}{2}}\right\rVert^{2} \\
   &\geq -C\Delta t\sum_{l=0}^{n-1}\left\lVert \widetilde{U}^{l+\frac{1}{2}}\right\rVert^{2}- \epsilon\Delta t\left\lVert\widetilde{U}^{n+\frac{1}{2}}\right\rVert^{2},
  \end{split}
\end{equation}
\begin{equation}\label{e12}
  \begin{split}
  \hspace{-3.5cm} T_{3}\leq \frac{(b^{*})^{2}+1}{2a_{*}^{2}}\Delta t \left\lVert \widetilde{U}^{n+\frac{1}{2}}\right\rVert^{2},
  \end{split}
\end{equation}
and
\begin{equation}\label{e13}
  \begin{split}
    \left(\widetilde{\psi}_{x} f^{n+\frac{1}{2}}, P^{n+\frac{1}{2}}\right)\leq \frac{1}{2}\left\lVert \widetilde{\psi}_{x}f^{n+\frac{1}{2}}\right\rVert^{2}+\frac{1}{2}\left\lVert P^{n+\frac{1}{2}}\right\rVert^{2}.
  \end{split}
\end{equation}
Taking $h$, $\epsilon$ and $\Delta t$ sufficiently small such that
\begin{equation*}
  \begin{split}
    \frac{5-6C_{a}h}{6a^{*}}-\epsilon\Delta t- \frac{(b^{*})^{2}+1}{2a_{*}^{2}}\Delta t \textgreater 0,
  \end{split}
\end{equation*}
multiplying Eq.$(\ref{e7})$ by $\Delta t$, summing on $n$ from 0 to $N-1$ and combining Eqs.$(\ref{e7})-(\ref{e13}),$ we have that
\begin{equation*}
  \begin{split}
    \frac{3}{4}\left\lVert P^{N}\right\rVert^{2} +C\sum_{n=0}^{N}\Delta t\left\lVert \widetilde{U}^{n+\frac{1}{2}}\right\rVert^{2}
    &\leq\left(\frac{5}{12\epsilon}+\frac{1}{2}\right)\Delta t\sum_{n=0}^{N}\left\lVert P^{n}\right\rVert^{2}+\epsilon\Delta t \sum_{n=0}^{N-1}\left\lVert \partial_{t}P^{n+1}\right\rVert^{2}\\
    &+\frac{\Delta t}{2}\sum_{n=0}^{N}\left\lVert \widetilde{\psi}_{x}f^{n}\right\rVert^{2}+\frac{C}{2}\Delta t\sum_{n=1}^{N-1}\Delta t\sum_{l=0}^{n-1}\left\lVert \widetilde{U}^{l+\frac{1}{2}}\right\rVert^{2}+\frac{4}{3}\left\lVert P^{0}\right\rVert^{2}.  
  \end{split}
\end{equation*}
By using Gronwall's inequality, we obtain
\begin{equation}\label{e14}
  \begin{split}
    \frac{3}{4}\left\lVert P^{N}\right\rVert^{2} +C\sum_{n=0}^{N}\Delta t\left\lVert \widetilde{U}^{n+\frac{1}{2}}\right\rVert^{2}
    \leq \frac{4}{3}\left\lVert P^{0}\right\rVert^{2}+\frac{\Delta t}{2}\sum_{n=0}^{N}\left\lVert \widetilde{\psi}_{x}f^{n}\right\rVert^{2}+\epsilon\Delta t \sum_{n=0}^{N-1}\left\lVert \partial_{t}P^{n+1}\right\rVert^{2}.  
  \end{split}
\end{equation}
Now, we multilpy the Eq.$(\ref{e4})$ by $2\Delta t\partial_{t}p^{n+1},$ sum on n from 0 to N-1, and get the following equation,
\begin{equation}\label{e15}
  \begin{split}
    2\Delta t\sum_{n=0}^{N-1}\left(\partial_{t} \widehat{\psi}_{x} P^{n+1}, \partial_{t}P^{n+1}\right)+2\Delta t\sum_{n=0}^{N-1}\left(\delta_{x}U^{n+\frac{1}{2}}, \partial_{t}P^{n+1}\right)=2\Delta t\sum_{n=0}^{N-1}\left(\widetilde{\psi}_{x} f^{n+\frac{1}{2}}, \partial_{t}P^{n+1}\right).
  \end{split}
\end{equation}
It is obvious that
\begin{equation}\label{e16}
  \begin{split}
    \left(\partial_{t} \widehat{\psi}_{x} P^{n+1}, \partial_{t}P^{n+1}\right)\geq \frac{3}{4}\left\lVert \partial_{t}P^{n+1}\right\rVert^{2}.
  \end{split}
\end{equation}
Next, we consider the second term on the left hand side of Eq.$(\ref{e15})$. 
\begin{equation*}
  \begin{split}
    \left(\delta_{x}U^{n+\frac{1}{2}}, \partial_{t}P^{n+1}\right)&=\left(U^{n+\frac{1}{2}},\psi_{x}\left(\frac{\partial_{t}\widetilde{U}}{a}^{n+1}\right)\right)\\
    &=\left(\frac{1}{2\Delta t}\left(\widetilde{U}^{n+1},\psi_{x}\left(\frac{\widetilde{U}}{a}\right)^{n+1}\right)-\frac{1}{2\Delta t}\left(\widetilde{U}^{n+1},\psi_{x}\left(\frac{\widetilde{U}}{a}\right)^{n}\right)\right.\\
    &\left.+\frac{1}{2\Delta t}\left(\widetilde{U}^{n},\psi_{x}\left(\frac{\widetilde{U}}{a}\right)^{n+1}\right)-\frac{1}{2\Delta t}\left(\widetilde{U}^{n},\psi_{x}\left(\frac{\widetilde{U}}{a}\right)^{n}\right)\right)\\
    &+\left(\Delta t\sum_{l=0}^{n}\left(\frac{b}{a}\widetilde{U}\right)^{l+\frac{1}{2}},\psi_{x}\left(\frac{\partial_{t}\widetilde{U}^{n+1}}{a}\right)\right)
    -\frac{1}{2}\left(\Delta t\left(\frac{b}{a}\widetilde{U}\right)^{n+\frac{1}{2}},\psi_{x}\left(\frac{\partial_{t}\widetilde{U}^{n+1}}{a}\right)\right)\\
    &=I_{1}+I_{2}+I_{3}.
  \end{split}
\end{equation*}
Therefore,
\begin{equation}\label{e17}
  \begin{split}
    2\Delta t\sum_{n=0}^{N-1}\left(\delta_{x}U^{n+\frac{1}{2}}, \partial_{t}P^{n+1}\right)=2\Delta t\sum_{n=0}^{N-1}I_{1}+2\Delta t\sum_{n=0}^{N-1}I_{2}+2\Delta t\sum_{n=0}^{N-1}I_{3}.
  \end{split}
\end{equation}
By the Lemma $\ref{addlem2}$, it is obvious that
\begin{equation}\label{e18}
  \begin{split}
    2\Delta t\sum_{n=0}^{N-1}I_{1}\geq &\left[\frac{5}{6a^{*}}-\left(\frac{C_{a}}{a^{*}}+\frac{D_{a}}{24}\right)h\right]\left\lVert \widetilde{U}^{N}\right\rVert^{2}-\left(\frac{1}{a_{*}}+\frac{D_{a}h}{24}\right)\left\lVert \widetilde{U}^{0}\right\rVert^{2}\\   
    &-\frac{8\left(a^{*}\right)^{2}\epsilon}{33}\sum_{n=0}^{N-1}\Delta t\left\lVert \partial_{t}P^{n+1}\right\rVert^{2}-\frac{D_{a}^{2}}{24\epsilon}\sum_{n=0}^{N}\Delta t\left\lVert \widetilde{U}^{n}\right\rVert^{2}.  
  \end{split}
\end{equation}
Then, we consider the last two terms on the right hand side of Eq.$(\ref{e17})$.
\begin{equation}\label{e19}
  \begin{split}
    2\Delta t\sum_{n=0}^{N-1}I_{2}&=2\Delta t\sum_{n=0}^{N-1}\Delta t\sum_{l=0}^{n}\left(\left(\frac{b}{a}\widetilde{U}\right)^{l+\frac{1}{2}},\psi_{x}\left(\frac{\partial_{t}\widetilde{U}^{n+1}}{a}\right)\right)\\
    &=2\Delta t\sum_{n=1}^{N}\left(\psi_{x}\left(\frac{\widetilde{U}}{a}\right)^{N}-\psi_{x}\left(\frac{\widetilde{U}}{a}\right)^{n-1},\left(\frac{b}{a}\widetilde{U}\right)^{n-\frac{1}{2}}\right)\\
    &\leq 2\Delta t\sum_{n=1}^{N}\epsilon\left(\left\lVert \psi_{x}\left(\frac{\widetilde{U}}{a}\right)^{N}-\psi_{x}\left(\frac{\widetilde{U}}{a}\right)^{n-1}\right\rVert^{2} \right)+\frac{2\Delta t}{4\epsilon}\sum_{n=1}^{N}\left\lVert \left(\frac{b}{a}\widetilde{U}\right)^{n-\frac{1}{2}}\right\rVert^{2}\\
    &\leq 2T\epsilon \left\lVert \psi_{x}\left(\frac{\widetilde{U}}{a}\right)^{N}\right\rVert^{2}+2\Delta t \epsilon\sum_{n=1}^{N} \left\lVert \psi_{x}\left(\frac{\widetilde{U}}{a}\right)^{n-1}\right\rVert^{2}+\frac{\Delta t}{2\epsilon}\sum_{n=1}^{N}\left\lVert \left(\frac{b}{a}\widetilde{U}\right)^{n-\frac{1}{2}}\right\rVert^{2}\\
    &\leq \frac{2T\epsilon}{a_{*}^{2}}\left\lVert \widetilde{U}^{N}\right\rVert^{2}+\frac{2\Delta t \epsilon}{a_{*}^{2}}\sum_{n=1}^{N}\left\lVert \widetilde{U}^{n-1}\right\rVert^{2}+\frac{\Delta t(b^{*})^{2}}{2\epsilon a_{*}^{2}}\sum_{n=1}^{N}\left\lVert \left(\widetilde{U}\right)^{n-\frac{1}{2}}\right\rVert^{2} \\
    &\leq  \frac{2T\epsilon}{a_{*}^{2}}\left\lVert \widetilde{U}^{N}\right\rVert^{2}+C\Delta t\sum_{n=0}^{N}\left\lVert \widetilde{U}^{n}\right\rVert^{2}. 
  \end{split}
\end{equation}
\begin{equation}\label{add2}
  \begin{split}
    \hspace{-5cm}2\Delta t\sum_{n=0}^{N-1}I_{3}&=\Delta t\sum_{n=0}^{N-1}\left(\left(\frac{b}{a}\widetilde{U}\right)^{n+\frac{1}{2}},\psi_{x}\left(\frac{\widetilde{U}^{n+1}}{a}\right)-\psi_{x}\left(\frac{\widetilde{U}^{n}}{a}\right)\right)\\
    &\leq C\Delta t\sum_{n=0}^{N-1}\left(\left\lVert \widetilde{U}^{n+1}\right\rVert^{2}+\left\lVert \widetilde{U}^{n}\right\rVert^{2}\right)\\
    &\leq C\Delta t\sum_{n=0}^{N}\left\lVert \widetilde{U}^{n}\right\rVert^{2}.
  \end{split}
\end{equation}
Using the Cauchy Schwarz inequality, we have
\begin{equation}\label{e20}
  \begin{split}
    2\Delta t\sum_{n=0}^{N-1}\left(\widetilde{\psi}_{x} f^{n+\frac{1}{2}}, \partial_{t}P^{n+1}\right)\leq \frac{1}{\epsilon}\sum_{n=0}^{N}\Delta t\left\lVert \widetilde{\psi}_{x}f^{n}\right\rVert^{2}+\epsilon\sum_{n=0}^{N}\Delta t \left\lVert \partial_{t}P^{n+1}\right\rVert^{2}. 
  \end{split}
\end{equation}
Combining Eqs.$(\ref{e15})-(\ref{e20}),$ we obtain that
\begin{equation}\label{e21}
  \begin{split}
    &\left(\frac{3}{2}-\frac{8\left(a*\right)^{2}\epsilon}{33}-\epsilon\right)\sum_{n=0}^{N-1}\Delta t \left\lVert \partial_{t}P^{n+1}\right\rVert^{2}+\left(\frac{5}{6a^{*}}-\left(\frac{C_{a}}{a^{*}}+\frac{D_{a}}{24}\right)h-\frac{2T\epsilon}{a_{*}^{2}}\right)\left\lVert \widetilde{U}^{N}\right\rVert^{2} \\
    &\leq \left(\frac{1}{a_{*}}+\frac{D_{a}h}{24}\right)\left\lVert \widetilde{U}^{0}\right\rVert^{2}+\left(\frac{D_{a}^{2}}{24\epsilon}+C\right)\sum_{n=0}^{N}\Delta t\left\lVert \widetilde{U}^{n}\right\rVert^{2} +\frac{1}{\epsilon}\sum_{n=0}^{N}\Delta t\left\lVert \widetilde{\psi}_{x}f^{n}\right\rVert^{2}.   
  \end{split}
\end{equation}
Taking $h$, $\epsilon$ and $\epsilon_{0}$ sufficiently small such that
\begin{equation*}
  \begin{split}
\frac{3}{2}-\frac{8\left(a*\right)^{2}\epsilon}{33}-\epsilon-\epsilon_{0}>0,
\end{split}
\end{equation*}
\begin{equation*}
  \begin{split}
    \frac{5}{6a^{*}}-\left(\frac{C_{a}}{a^{*}}+\frac{D_{a}}{24}\right)h-\frac{2T\epsilon}{a_{*}^{2}}>0,
\end{split}
\end{equation*}
using the Gronwall's inequality and combining Eq.$(\ref{e14})$ and Eq.$(\ref{e21})$, we can easily get
\begin{equation*}
  \begin{split}
    \left\lVert P^{N}\right\rVert^{2}+\left\lVert \widetilde{U}^{N}\right\rVert^{2}\leq C\left(\left\lVert P^{0}\right\rVert^{2}+\left\lVert \widetilde{U}^{0}\right\rVert^{2}+\sum_{n=0}^{N}\Delta t\left\lVert \widetilde{\psi}_{x}f^{n}\right\rVert^{2} \right).  
\end{split}
\end{equation*}
Therefore, we get the following stability theorem.
\begin{thm}\label{thm1}
  Suppose $h$ and $\Delta t$ are sufficiently small. Let $P$ and $\widetilde{U}$ be the solutions to CBCFD scheme. There exists a positive constant $C$ such that the following inequality holds
  \begin{equation}\label{e22}
    \begin{split}
      \left\lVert P^{N}\right\rVert^{2}+\left\lVert \widetilde{U}^{N}\right\rVert^{2}\leq C\left(\left\lVert P^{0}\right\rVert^{2}+\left\lVert \widetilde{U}^{0}\right\rVert^{2}+\sum_{n=0}^{N}\Delta t\left\lVert \widetilde{\psi}_{x}f^{n}\right\rVert^{2} \right).  
  \end{split}
  \end{equation}
\end{thm}
\subsection{The error estimate}
The convergence of the CBCFD scheme for the one-dimensional non-Fickian model will be discussed in this section. Suppose $b(x,t),p, u\in C^{3}(0,T;C^{5}(\Omega))$ and $a(x)\in C^{5}(\Omega).$ By the Eq.$(\ref{e3})$ and Taylor's expansion, the following equations can be easily obtained.
\begin{equation}\label{e23}
  \begin{split}
  \partial_{t} \widehat{\psi}_{x} p_{i}^{n+1}+\delta_{x} u_{i}^{n+\frac{1}{2}}=\widetilde{\psi}_{x} f_{i}^{n+\frac{1}{2}}+O\left(h^{4}+\Delta t^{2}\right)
  \end{split}
  \end{equation}
  \begin{equation}\label{e24}
    \begin{split}
    \delta_{x} p_{i+1 / 2}^{n+1}=-\psi_{x}\left(\frac{\widetilde{u}}{a}\right)_{i+1 / 2}^{n+1}+O\left(h^{4}\right)
    \end{split}
  \end{equation}
  \begin{equation}\label{e25}
    \begin{split}
    \partial_{t}\delta_{x} p_{i+1 / 2}^{n+1}=-\partial_{t}\psi_{x}\left(\frac{\widetilde{u}}{a}\right)_{i+1 / 2}^{n+1}+O\left(h^{4}\right)
    \end{split}
  \end{equation}
Set
\begin{equation*}
    \begin{split}
    \eta&=P-p,\\
    \widetilde{\xi}&=\widetilde{U}-\widetilde{u},\\
    \xi&=U-u=U-\widetilde{U}+\widetilde{\xi}+\widetilde{u}-u.
    \end{split}
\end{equation*}
Subtracting Eq.$(\ref{e23})$ from Eq.$(\ref{e4})$ and Eq.$(\ref{e24})$ from Eq.$(\ref{e5})$, we get the error equations
\begin{equation}\label{e26}
  \begin{split}
  \partial_{t} \widehat{\psi}_{x} \eta_{i}^{n+1}+\delta_{x} \xi_{i}^{n+\frac{1}{2}}=R_{1,i}^{n+1},
  \end{split}
  \end{equation}
  \begin{equation}\label{e27}
    \begin{split}
    \delta_{x} \eta_{i+1 / 2}^{n+1}=-\psi_{x}\left(\frac{\widetilde{\xi}}{a}\right)_{i+1 / 2}^{n+1}+R_{2,i+\frac{1}{2}}^{n+1},
    \end{split}
  \end{equation}
  \begin{equation}\label{add1}
    \begin{split}
    \partial_{t}\delta_{x} \eta_{i+1 / 2}^{n+1}=-\partial_{t}\psi_{x}\left(\frac{\widetilde{\xi}}{a}\right)_{i+1 / 2}^{n+1}+R_{3,i+\frac{1}{2}}^{n+1},
    \end{split}
  \end{equation}
  where $R_{1,i}^{n+1}=O(h^{4}+\Delta t^{2}),$ $R_{2,i+\frac{1}{2}}^{n+1}=O(h^{4})$ and $R_{3,i+\frac{1}{2}}^{n+1}=O(h^{4})$ .

  Multiplying Eq.$(\ref{e26})$ by $h_{i}\eta_{i}^{n+\frac{1}{2}}$ and summing on $i$, we have that
  \begin{equation}\label{e28}
    \begin{split}
    &\left(\partial_{t} \widehat{\psi}_{x} \eta^{n+1}, \eta^{n+\frac{1}{2}}\right)+\left(\xi^{n+\frac{1}{2}}, \psi_{x}\left(\frac{\widetilde{\xi}}{a}\right)^{n+\frac{1}{2}}\right) \\
    &=\left(R_{1}^{n+1}, \eta^{n+\frac{1}{2}}\right)+\left(\xi^{n+\frac{1}{2}}, R_{2}^{n+1}\right).
    \end{split}
    \end{equation}
Similar to the stability analysis, we can get the following estimates,
\begin{equation}\label{e29}
  \begin{split}
  &\left(\partial_{t} \widehat{\psi}_{x} \eta^{n+1}, \eta^{n+\frac{1}{2}}\right)=\frac{1}{\Delta t}\left(\left(\widehat{\psi}_{x}\eta^{n+1},\eta^{n+1}\right)-\left(\widehat{\psi}_{x}\eta^{n},\eta^{n}\right)\right)-\left(\widehat{\psi}_{x}\eta^{n+\frac{1}{2}},\partial_{t}\eta^{n+1}\right),
  \end{split}
\end{equation}
  \begin{equation}\label{e30}
    \begin{split}
      \left(R_{1}^{n+1}, \eta^{n+\frac{1}{2}}\right)\leq C\left(\Delta t^{4}+h^{8}+\left\lVert \eta^{n+\frac{1}{2}}\right\rVert^{2} \right),
    \end{split}
  \end{equation}
  \begin{equation}\label{e31}
    \begin{split}
      \left(R_{2}^{n+1}, \xi^{n+\frac{1}{2}}\right)\leq Ch^{8}+\epsilon\left\lVert \xi^{n+\frac{1}{2}}\right\rVert^{2} .
    \end{split}
  \end{equation}
Next, we consider the estimate of the second term on the left hand side of Eq.$(\ref{e28})$.
\begin{equation}\label{e32}
  \begin{split}
    \left(\xi^{n+\frac{1}{2}}, \psi_{x}\left(\frac{\widetilde{\xi}}{a}\right)^{n+\frac{1}{2}}\right)
    &=\left(U^{n+\frac{1}{2}}-\widetilde{U}^{n+\frac{1}{2}}+\widetilde{\xi}^{n+\frac{1}{2}}+\widetilde{u}^{n+\frac{1}{2}}-u^{n+\frac{1}{2}},\psi_{x}\left(\frac{\widetilde{\xi}}{a}\right)^{n+\frac{1}{2}}\right)\\
    &=\left(\frac{\Delta t}{2}\sum_{l=0}^{n}\left(\frac{b}{a}\widetilde{U}\right)^{l+\frac{1}{2}}-\frac{1}{2}\int_{0}^{t_{n+1}}\frac{b}{a}\widetilde{u}\left(x,s\right)ds,\psi_{x}\left(\frac{\widetilde{\xi}}{a}\right)^{n+\frac{1}{2}}\right)\\
    &\quad+\left(\frac{\Delta t}{2}\sum_{l=0}^{n-1}\left(\frac{b}{a}\widetilde{U}\right)^{l+\frac{1}{2}}-\frac{1}{2}\int_{0}^{t_{n}}\frac{b}{a}\widetilde{u}\left(x,s\right)ds,\psi_{x}\left(\frac{\widetilde{\xi}}{a}\right)^{n+\frac{1}{2}}\right)\\
    &\quad+\left(\widetilde{\xi}^{n+\frac{1}{2}}, \psi_{x}\left(\frac{\widetilde{\xi}}{a}\right)^{n+\frac{1}{2}}\right)\\
    &=S_{1}+S_{2}+S_{3}
  \end{split}
\end{equation}
By the straightforward calculation, we obtain that
\begin{equation}\label{add3}
  \begin{split}
    S_{1}&=\frac{1}{2}\left(\Delta t\sum_{l=0}^{n}\left(\frac{b}{a}\left(\widetilde{U}-\widetilde{u}\right)\right)^{l+\frac{1}{2}}+\Delta t\sum_{l=0}^{n}\left(\frac{b}{a}\widetilde{u}\right)^{l+\frac{1}{2}}-\int_{0}^{t_{n+1}}\frac{b}{a}\widetilde{u}\left(x,s\right)ds,\psi_{x}\left(\frac{\widetilde{\xi}}{a}\right)^{n+\frac{1}{2}}\right)\\
    &=\frac{1}{2}\left(\Delta t\sum_{l=0}^{n}\left(\frac{b}{a}\widetilde{\xi}\right)^{l+\frac{1}{2}}+\sum_{l=0}^{n}\int_{t_{l}}^{t_{l+1}}\left(\frac{b}{a}\widetilde{u}\right)^{l+\frac{1}{2}}dt-\int_{0}^{t_{n+1}}\frac{b}{a}\widetilde{u}\left(x,s\right)ds,\psi_{x}\left(\frac{\widetilde{\xi}}{a}\right)^{n+\frac{1}{2}}\right)\\
    &\leq C\Delta t\sum_{l=0}^{n}\left\lVert \widetilde{\xi}^{l+\frac{1}{2}}\right\rVert^{2}+O(\Delta t^{4})+\epsilon\left\lVert \widetilde{\xi}^{n+\frac{1}{2}}\right\rVert^{2}.  
  \end{split}
\end{equation}
Similarly,
\begin{equation}\label{add4}
  \begin{split}
    S_{2}\leq C\Delta t\sum_{l=0}^{n-1}\left\lVert \widetilde{\xi}^{l+\frac{1}{2}}\right\rVert^{2}+O(\Delta t^{4})+\epsilon\left\lVert \widetilde{\xi}^{n+\frac{1}{2}}\right\rVert^{2}.  
  \end{split}
\end{equation}
By using Lemma $\ref{addlem1},$ we have
\begin{equation*}
  \begin{split}
   S_{3}\geq \frac{5-6C_{a}h}{6a^{*}}\left\lVert \widetilde{\xi}^{n+\frac{1}{2}}\right\rVert^{2}. 
  \end{split}
\end{equation*}
Therefore, Eq.$(\ref{e32})$ can be estimated as
\begin{equation}\label{e33}
  \begin{split}
    \left(\xi^{n+\frac{1}{2}}, \psi_{x}\left(\frac{\widetilde{\xi}}{a}\right)^{n+\frac{1}{2}}\right)\geq \left(\frac{5-6C_{a}h}{6a^{*}}-2\epsilon\right)\left\lVert \widetilde{\xi}^{n+\frac{1}{2}}\right\rVert^{2}-C\Delta t\sum_{l=0}^{n}\left\lVert \widetilde{\xi}^{l+\frac{1}{2}}\right\rVert^{2}-O(\Delta t^{4}).
  \end{split}
\end{equation}
Multiplying Eq.$(\ref{e28})$ by $\Delta t$, combining the Eqs.$(\ref{e28})-(\ref{e31})$ and Eq.$(\ref{e33})$, and summing on n from $0$ to $N-1$, we obtain the following inequality,
\begin{equation}\label{e34}
  \begin{split}
   \frac{3}{4}\left\lVert \eta^{N}\right\rVert^{2}+&\left(\frac{5-6C_{a}h}{6a^{*}}-3\epsilon\right)\sum_{n=0}^{N-1}\Delta t\left\lVert \widetilde{\xi}^{n+\frac{1}{2}}\right\rVert^{2}\\
   \leq &\frac{4}{3}\left\lVert \eta^{0}\right\rVert^{2}+\epsilon\sum_{n=0}^{N-1}\Delta t\left\lVert \partial_{t}\eta^{n+1}\right\rVert^{2}
+C\sum_{n=0}^{N-1}\Delta t\left\lVert \eta^{n+\frac{1}{2}}\right\rVert^{2}\\ &+C\Delta t\sum_{n=0}^{N-1}\Delta t\sum_{l=0}^{n}\left\lVert \widetilde{\xi}^{l+\frac{1}{2}}\right\rVert^{2}
   +O(\Delta t^{4}+h^{8}).    
  \end{split}
\end{equation}
Taking $h$ and $\epsilon$ sufficiently small such that
\begin{equation*}
  \begin{split}
    \frac{5-6C_{a}h}{6a^{*}}-3\epsilon>0,  
  \end{split}
\end{equation*}
and using the Gronwall's inequality, we can obtain
\begin{equation}\label{e35}
  \begin{split}
   \frac{3}{4}\left\lVert \eta^{N}\right\rVert^{2}+&\left(\frac{5-6C_{a}h}{6a^{*}}-3\epsilon\right)\sum_{n=0}^{N-1}\Delta t\left\lVert \widetilde{\xi}^{n+\frac{1}{2}}\right\rVert^{2}\\
   \leq &\frac{4}{3}\left\lVert \eta^{0}\right\rVert^{2}+\epsilon\sum_{n=0}^{N-1}\Delta t\left\lVert \partial_{t}\eta^{n+1}\right\rVert^{2}+O(\Delta t^{4}+h^{8}).    
  \end{split}
\end{equation}
Multilpying Eq.$(\ref{e26})$ by $2\Delta th_{i}\partial_{t}\eta_{i}^{n+1},$ summing on $i$ and $n$ , we have that
  \begin{equation}\label{e37}
    \begin{split}
    &2\Delta t\sum_{n=0}^{N-1}\left(\partial_{t} \widehat{\psi}_{x} \eta^{n+1}, \partial_{t}\eta^{n+1}\right)+2\Delta t\sum_{n=0}^{N-1}\left(\xi^{n+\frac{1}{2}}, \partial_{t}\psi_{x}\left(\frac{\widetilde{\xi}}{a}\right)^{n+1}\right) \\
    &=2\Delta t\sum_{n=0}^{N-1}\left(R_{1}^{n+1}, \partial_{t}\eta^{n+1}\right)+2\Delta t\sum_{n=0}^{N-1}\left(\xi^{n+\frac{1}{2}}, R_{3}^{n+1}\right).
    \end{split}
    \end{equation}
The following estimates can be easily got.
\begin{equation}\label{e38}
  \begin{split}
  \left(\partial_{t} \widehat{\psi}_{x} \eta^{n+1}, \partial_{t}\eta^{n+1}\right)\geq \frac{3}{4}\left\lVert \partial_{t}\eta^{n+1}\right\rVert^{2}.
  \end{split}
  \end{equation}
  \begin{equation}\label{e39}
    \begin{split}
      \left(R_{1}^{n+1}, \partial_{t}\eta^{n+1}\right)\leq O\left(\Delta t^{4}+h^{8}\right)+\widetilde{\epsilon}\left\lVert \partial_{t}\eta^{n+1}\right\rVert^{2}.
    \end{split}
    \end{equation}
    \begin{equation}\label{e40}
      \begin{split}
        \left(\xi^{n+\frac{1}{2}}, R_{3}^{n+1}\right)\leq O(h^{8})+C\left\lVert \xi^{n+\frac{1}{2}}\right\rVert^{2}.
      \end{split}
      \end{equation}
Similar to the estimates of the Eq.$(\ref{e18})$, $(\ref{add3})$ and $(\ref{add4})$, we obtain the following estimate by using the Lemma $\ref{addlem3}.$  
\begin{equation}\label{e44}
  \begin{split}
&2\Delta t\sum_{n=0}^{N-1}\left(\xi^{n+\frac{1}{2}}, \partial_{t}\psi_{x}\left(\frac{\widetilde{\xi}}{a}\right)^{n+1}\right)\\
&=2\Delta t\sum_{n=0}^{N-1}\left(\widetilde{\xi}^{n+\frac{1}{2}}, \partial_{t}\psi_{x}\left(\frac{\widetilde{\xi}}{a}\right)^{n+1}\right)+2\Delta t\sum_{n=0}^{N-1}\left(U^{n+\frac{1}{2}}-\widetilde{U}^{n+\frac{1}{2}}, \partial_{t}\psi_{x}\left(\frac{\widetilde{\xi}}{a}\right)^{n+1}\right)\\
&\quad+2\Delta t\sum_{n=0}^{N-1}\left(\widetilde{u}^{n+\frac{1}{2}}-u^{n+\frac{1}{2}}, \partial_{t}\psi_{x}\left(\frac{\widetilde{\xi}}{a}\right)^{n+1}\right)\\
&\geq \left[\frac{5}{6a^{*}}-\left(\frac{C_{a}}{a^{*}}+\frac{D_{a}}{24}\right)h\right]\left\lVert \widetilde{\xi}^{N}\right\rVert^{2}-\left(\frac{1}{a_{*}}+\frac{D_{a}h}{24}\right)\left\lVert \widetilde{\xi}^{0}\right\rVert^{2}-\frac{8\left(a^{*}\right)^{2}\epsilon_{0}}{33}\sum_{n=0}^{N-1}\Delta t\left\lVert \partial_{t}\eta^{n+1}\right\rVert^{2}\\
&\quad-\frac{D_{a}^{2}}{24\epsilon}\sum_{n=0}^{N}\Delta t\left\lVert \widetilde{\xi}^{n}\right\rVert^{2}-C\Delta t\sum_{n=0}^{N-1}\Delta t\sum_{l=0}^{n}\left\lVert \widetilde{\xi}^{l+\frac{1}{2}}\right\rVert^{2}-O(\Delta t^{4})-2\epsilon_{1}\sum_{n=0}^{N-1}\Delta t\left\lVert \partial_{t}\psi_{x}\left(\frac{\widetilde{\xi}}{a}\right)^{n+1}\right\rVert^{2}\\
&\geq \left[\frac{5}{6a^{*}}-\left(\frac{C_{a}}{a^{*}}+\frac{D_{a}}{24}\right)h\right]\left\lVert \widetilde{\xi}^{N}\right\rVert^{2}-\left(\frac{1}{a_{*}}+\frac{D_{a}h}{24}\right)\left\lVert \widetilde{\xi}^{0}\right\rVert^{2}-\left(\frac{8\left(a^{*}\right)^{2}\epsilon_{0}}{33}+\frac{8\epsilon_{1}}{h^{2}}\right)\sum_{n=0}^{N-1}\Delta t\left\lVert \partial_{t}\eta^{n+1}\right\rVert^{2}\\
&\quad-\frac{D_{a}^{2}}{24\epsilon}\sum_{n=0}^{N}\Delta t\left\lVert \widetilde{\xi}^{n}\right\rVert^{2}-C\Delta t\sum_{n=0}^{N-1}\Delta t\sum_{l=0}^{n}\left\lVert \widetilde{\xi}^{l+\frac{1}{2}}\right\rVert^{2}-O(\Delta t^{4}).
\end{split}
\end{equation}
Combining the Eqs.$(\ref{e37})-(\ref{e44}),$ we have that 
\begin{equation}\label{e45}
  \begin{split}
&\frac{3}{2}\Delta t\sum_{n=0}^{N-1}\left\lVert \partial_{t}\eta^{n+1}\right\rVert^{2}+\left[\frac{5}{6a^{*}}-\left(\frac{C_{a}}{a^{*}}+\frac{D_{a}}{24}\right)h\right]\left\lVert \widetilde{\xi}^{N}\right\rVert^{2}\\
&\leq \left(\frac{1}{a_{*}}+\frac{D_{a}h}{24}\right)\left\lVert \widetilde{\xi}^{0}\right\rVert^{2}+\left(\frac{8\left(a^{*}\right)^{2}\epsilon_{0}}{33}+\frac{8\epsilon_{1}}{h^{2}}+2\widetilde{\epsilon}\right)\sum_{n=0}^{N-1}\Delta t\left\lVert \partial_{t}\eta^{n+1}\right\rVert^{2}+\left(\frac{D_{a}^{2}}{24\epsilon}+C\right)\sum_{n=0}^{N}\Delta t\left\lVert \widetilde{\xi}^{n}\right\rVert^{2}\\
&+C\Delta t\sum_{n=0}^{N-1}\Delta t\sum_{l=0}^{n}\left\lVert \widetilde{\xi}^{l+\frac{1}{2}}\right\rVert^{2}+O(\Delta t^{4}+h^{8}).
\end{split}
\end{equation}
Taking $\epsilon_{0},\epsilon_{1},\widetilde{\epsilon}$ and $\epsilon$ sufficiently small such that
\begin{equation*}
  \begin{split}
    \frac{3}{2}-\frac{8\left(a^{*}\right)^{2}\epsilon_{0}}{33}+\frac{8\epsilon_{1}}{h^{2}}+2\widetilde{\epsilon}-\epsilon>0,
  \end{split}
\end{equation*}
 combining the Eq.$(\ref{e35})$ and Eq.$(\ref{e45})$ and using Gronwall's inequality, we obtain that
 \begin{equation}\label{e46}
  \begin{split}
    \left\lVert \eta^{N}\right\rVert^{2}+\left\lVert \widetilde{\xi}^{N}\right\rVert^{2}\leq O\left(\Delta t^{4}+h^{8}\right).
\end{split}
\end{equation}
Then, we get the convergence conclusion.
\begin{thm}\label{thm2}
  Suppose $h$ and $\Delta t$ are sufficiently small and $p, u\in C^{3}(0,T;C^{5}(\Omega))$. Let $P$ and $\widetilde{U}$ be the numerical solutions to the CBCFD scheme. There exists a positive constant $C$ such that the following inequality holds,
  \begin{equation}\label{e22}
    \begin{split}
      \left\lVert P^{N}-p^{N}\right\rVert+\left\lVert \widetilde{U}^{N}-\widetilde{u}^{N}\right\rVert\leq C\left(\Delta t^2+h^{4} \right).  
  \end{split}
  \end{equation}
\end{thm}
\section{Two-dimensional parabolic problem}
Now, we consider the two-dimensional non-Fickian flow model with variable coefficients. The problem can be expressed as
\begin{equation}\label{e48}
  \left\{\begin{array}{ll} p_{t}+\nabla\cdot \boldsymbol{u}=f(x,y,t),  &(x, y) \in \Omega, t \in J, \\
  \boldsymbol{u}=-\left(\boldsymbol{A}\nabla p+\int_{0}^{t}\boldsymbol{B}\left(s\right)\nabla p\left(s\right)ds\right), & (x, y) \in \Omega, t \in J,\\
  p|_{t=0}=p_{0}(x,y),  &(x,y) \in \Omega,
\end{array}\right.
\end{equation}
where $\Omega=(0,L_{1})\times(0,L_{2})$, J=(0,T], and the boundary conditions are periodic. Besides, the definitions of $\boldsymbol{A}$ and $\boldsymbol{B}$ are as follows
\begin{equation*}
  \begin{split}
 \boldsymbol{A}&=diag\left(a^{x}(x,y),a^{y}(x,y)\right),\\
 \boldsymbol{B}&=diag\left(b^{x}(x,y,t),b^{y}(x,y,t)\right).
  \end{split}
\end{equation*}
 We suppose that there exist four positive constants $a_{*},a^{*},b_{*}$ and $b^{*},$ such that
 $$0\textless a_{*}\leq a^{x},a^{y}\leq a^{*},\quad 0\textless b_{*}\leq b^{x},b^{y}\leq b^{*}. $$
 Besides, let $a^{x},$ $a^{y},$ $b^{x},$ $b^{y},$ and their partial derivatives with respect to $x$ and $y$ are bounded smooth functions.

Set $\widetilde{u}=-\boldsymbol{A}\nabla p$. The Eq.$(\ref{e48})$ can be recast into the following formulation:
\begin{equation}\label{e49}
  \left\{\begin{array}{ll} p_{t}+\nabla\cdot \boldsymbol{u}=f(x,t),  &(x, y) \in \Omega, t \in J, \\
  \boldsymbol{u}=\widetilde{\boldsymbol{u}}+\int_{0}^{t}\boldsymbol{B}(s)\boldsymbol{A}^{-1}(s)\widetilde{\boldsymbol{u}}ds & (x, y) \in \Omega, t \in J,\\
  p|_{t=0}=p_{0}(x,y),  &(x,y) \in \Omega.
\end{array}\right.
\end{equation}
\subsection{The CBCFD scheme of two dimensional problem}
We define $\Omega_{i,j}=[x_{i-1/2},x_{i+1/2}]\times [y_{j-1/2},y_{j+1/2}],$ $i=1, 2,\cdots, N_{1}$, $j=1, 2,\cdots, N_{2}$ as the uniform partitions of $\Omega,$ where $x_{1/2}=0,$ $x_{N_{1}+1/2}=L_{1},$ $y_{1/2}=0,$ and $y_{N_{2}+1/2}=L_{2}.$ Let $h=L_{1}/N_{1}$ and $k=L_{2}/N_{2}.$

Set 
$$x_{i}=\frac{x_{i-1/2}+x_{i+1/2}}{2},$$
and
$$y_{j}=\frac{y_{j-1/2}+y_{j+1/2}}{2}.$$

Let $f$ and $g$ be any functions, we define the discrete inner operators as follows:
\begin{equation*}
  \begin{aligned}
  &(f, g)=\sum_{i=1}^{N_{1}} \sum_{j=1}^{N_{2}}hkf_{i j}g_{i j},\\
  &(f, g)_{x}=\sum_{i=0}^{N_{1}-1} \sum_{j=1}^{N_{2}}hkf_{i+\frac{1}{2}, j} g_{i+\frac{1}{2}, j},\\
  &(f, g)_{y}=\sum_{i=1}^{N_{1}} \sum_{j=0}^{N_{2}-1} hkf_{i, j+\frac{1}{2}} g_{i, j+\frac{1}{2}}.
  \end{aligned}
  \end{equation*}
Define
\begin{equation*}
  \begin{gathered}
  \delta_{x} f_{i+\frac{1}{2}, j}=\frac{f_{i+1, j}-f_{i j}}{h}, \quad \delta_{y} f_{i, j+\frac{1}{2}}=\frac{f_{i, j+1}-f_{i j}}{k}, \\
  \delta_{x} f_{i j}=\frac{f_{i+\frac{1}{2}, j}-f_{i-\frac{1}{2}, j}}{h}, \quad \delta_{y} f_{i j}=\frac{f_{i, j+\frac{1}{2}}-f_{i, j-\frac{1}{2}}}{k}, \\
  \psi_{x}=\left(1+\frac{h^{2}}{24} \delta_{x}^{2}\right), \quad \psi_{y}=\left(1+\frac{k^{2}}{24} \delta_{y}^{2}\right).
  \end{gathered}
  \end{equation*}
The definitions of $\widetilde{\psi}_{x},$ $\widetilde{\psi}_{y},$ $\widehat{\psi}_{x},$ and $\widehat{\psi}_{y}$ can be easily obtained by the one-dimensional case. By using above definitions, we can establish the CBCFD scheme for the two-dimensional case.
\begin{equation}\label{e50}
  \begin{split}
  \frac{\widehat{\psi}_{x}\widehat{\psi}_{y}P^{n+1}_{i,j}-\widehat{\psi}_{x}\widehat{\psi}_{y}P^{n}_{i,j}}{\Delta t}+\widehat{\psi}_{y}\delta_{x}U^{x,n+\frac{1}{2}}_{i,j}+\widehat{\psi}_{x}\delta_{y}U^{y,n+\frac{1}{2}}_{i,j}=\widetilde{\psi}_{x}\widetilde{\psi}_{y}f^{n+\frac{1}{2}}_{i,j},
\end{split}
\end{equation}
\begin{equation}\label{e51}
  \begin{split}
   \delta_{x}P^{n+1}_{i+\frac{1}{2},j}=-\psi_{x}\left(\frac{\widetilde{U}^{x}}{a^{x}}\right)^{n+1}_{i+\frac{1}{2},j},
\end{split}
\end{equation}
\begin{equation}\label{e52}
  \begin{split}
   \delta_{y}P^{n+1}_{i,j+\frac{1}{2}}=-\psi_{y}\left(\frac{\widetilde{U}^{y}}{a^{y}}\right)^{n+1}_{i,j+\frac{1}{2}},
\end{split}
\end{equation}
\begin{equation}\label{e53}
  \begin{split}
   U^{x,n+1}_{i+\frac{1}{2},j}=\widetilde{U}^{x,n+1}_{i+\frac{1}{2},j}+\Delta t\sum_{l=0}^{n}\left(\frac{b^{x}}{a^{x}}\widetilde{U}^{x}\right)^{l+\frac{1}{2}}_{i+\frac{1}{2},j},
\end{split}
\end{equation}
\begin{equation}\label{e54}
  \begin{split}
   U^{y,n+1}_{i,j+\frac{1}{2}}=\widetilde{U}^{y,n+1}_{i,j+\frac{1}{2}}+\Delta t\sum_{l=0}^{n}\left(\frac{b^{y}}{a^{y}}\widetilde{U}^{y}\right)^{l+\frac{1}{2}}_{i,j+\frac{1}{2}}.
\end{split}
\end{equation}
Under the periodic boundary conditions, the difference operators $\widehat{\psi}_{x}$ (or $\widetilde{\psi}_{x}$) and $\widehat{\psi}_{y}$ (or $\widetilde{\psi}_{y}$) can be replaced by $\psi_{x}$ and $\psi_{y}$. Thus, Eq.$(\ref{e50})$ is equivalent to the following equation,
\begin{equation}\label{e55}
  \begin{split}
  \frac{\psi_{x}\psi_{y}P^{n+1}_{i,j}-\psi_{x}\psi_{y}P^{n}_{i,j}}{\Delta t}+\psi_{y}\delta_{x}U^{x,n+\frac{1}{2}}_{i,j}+\psi_{x}\delta_{y}U^{y,n+\frac{1}{2}}_{i,j}=\psi_{x}\psi_{y}f^{n+\frac{1}{2}}_{i,j}.
\end{split}
\end{equation}
Now, we give some lemmas which will be useful for the theoretical analysis. 
  \begin{lem}\label{lem2}
    Supposing that the boundary conditions are periodic and $P_{i,j}$ is the approximation to $p_{i,j}$, then we can obtain that
    \begin{equation*}
     \begin{split}
     \left(\psi_{x}\psi_{y}P,P\right)\geq \frac{49}{72}\left\lVert P\right\rVert^{2}. 
     \end{split}
   \end{equation*}
   
\begin{proof}
   
   \begin{equation*}
     \begin{split}
      \psi_{x}\psi_{y}P_{i,j}P_{i,j}=&\frac{\psi_{y}P_{i+1,j}+22\psi_{y}P_{i,j}+\psi_{y}P_{i-1,j}}{24}P_{i,j}\\
      =&\frac{1}{576}\left(P_{i+1,j+1}+22P_{i+1,j}+P_{i+1,j-1}+22P_{i,j+1}+484P_{i,j}+22P_{i,j-1}\right.\\
      &\left.+P_{i-1,j+1}+22P_{i-1,j}+P_{i-1,j-1}\right)P_{i,j}\\
      \geq & \frac{1}{576}\left(-\frac{1}{2}P_{i+1,j+1}^{2}-11P_{i+1,j}^{2}-\frac{1}{2}P_{i+1,j-1}^{2}-11P_{i,j+1}^{2}+438P_{i,j}^{2}-11P_{i,j-1}^{2}-\frac{1}{2}P_{i-1,j+1}^{2}-11P_{i-1,j}^{2}-\frac{1}{2}P_{i-1,j-1}^{2}\right).
     \end{split}
   \end{equation*}
   Therefore,
   \begin{equation*}
     \begin{split}
     \left(\psi_{x}\psi_{y}P,P\right)=\sum_{i=1}^{N_{1}}\sum_{j=1}^{N_{2}}hk\psi_{x}\psi_{y}P_{i,j}P_{i,j} \geq \frac{49}{72}\sum_{i=1}^{N_{1}}\sum_{j=1}^{N_{2}}hkP_{i,j}^{2}=\frac{49}{72}\left\lVert P \right\rVert^{2}. 
     \end{split}
   \end{equation*}
\end{proof}   
   
   \end{lem}
  \begin{lem}\label{lem3}\cite{2021High}
    Let $C_{a^{x}}=\max \left\{\left\|\frac{\partial a^{x}}{\partial x}\right\|_{\infty},\left\|\frac{\partial a^{x}}{\partial y}\right\|_{\infty}\right\}, C_{a^{y}}=\max \left\{\left\|\frac{\partial a^{y}}{\partial x}\right\|_{\infty},\left\|\frac{\partial a^{y}}{\partial y}\right\|_{\infty}\right\}$ and $h_{max}=max(h,k)$, we can obtain that
    \begin{equation*}
      \begin{split}
    \left(a^{x} V, \psi_{x} \psi_{y} V\right)_{x} & \geq \frac{49-6 C_{a^{x}} h_{max}}{72}\left(a^{x} V, V\right)_{x}, \\
    \left(a^{y} V, \psi_{x} \psi_{y} V\right)_{y} & \geq \frac{49-6 C_{a^{y}} h_{max}}{72}\left(a^{y} V, V\right)_{y},
      \end{split}
    \end{equation*}
  and 
  \begin{equation*}
    \begin{split}
    (\psi_{x}V,\psi_{x}V)_{x}\leq \left\lVert V\right\rVert_{x}^{2},\\
    (\psi_{y}V,\psi_{y}V)_{y}\leq \left\lVert V\right\rVert_{y}^{2}.
    \end{split}
  \end{equation*}
\end{lem}
\begin{lem}\label{lem7}\cite{2021High}
Assuming that $h,$ $k,$ and $\epsilon$ are sufficiently small, then there exists a positive constant $C$ such that
\begin{equation*}
  \begin{split}
2\Delta t\sum_{n=1}^{N-1}\left(\widetilde{U}^{s,n+\frac{1}{2}},\partial_{t}\psi_{x}\psi_{y}\left(\frac{\widetilde{U}^{s}}{a^{s}}\right)^{n+1}\right)\geq 
&C\left(\left\lVert \widetilde{U}^{s,N}\right\rVert_{s}^{2}-\left\lVert \widetilde{U}^{s,0}\right\rVert_{s}^{2}\right)-\epsilon\sum_{n=0}^{N-1}\Delta t\left\lVert \partial_{t}P^{n+1}\right\rVert^{2}\\
&-C\Delta t\sum_{n=0}^{N}\left\lVert \widetilde{U}^{s,n}\right\rVert_{s}^{2},
\end{split}
\end{equation*}
where $s=x$ or $y.$
\end{lem}

\subsection{The analysis of stability}
In this section, we will derive the stablity of the CBCFD scheme for two dimensional non-Fickian flow problem.\\ 
Multiplying Eq.$(\ref{e55})$ by $h_{i}k_{j}P_{i,j}^{n+\frac{1}{2}}$ and summing on $i$ and $j$, we obtain that
\begin{equation}\label{e56}
  \begin{split}
  \left(\partial_{t}\psi_{x}\psi_{y}P^{n+1},P^{n+\frac{1}{2}}\right)+\left(\psi_{y}\delta_{x}U^{x,n+\frac{1}{2}},P^{n+\frac{1}{2}}\right)+\left(\psi_{x}\delta_{y}U^{y,n+\frac{1}{2}},P^{n+\frac{1}{2}}\right)=\left(\psi_{x}\psi_{y}f^{n+\frac{1}{2}},P^{n+\frac{1}{2}}\right).
\end{split}
\end{equation}
Substituting Eqs.$(\ref{e51})-(\ref{e54})$ into the Eq.$(\ref{e55})$ and using the Lemma $\ref{lem1}$, we then have that
\begin{equation}\label{e57}
  \begin{split}
  &\left(\partial_{t}\psi_{x}\psi_{y}P^{n+1},P^{n+\frac{1}{2}}\right)+\left(\widetilde{U}^{x,n+\frac{1}{2}},\psi_{x}\psi_{y}\left(\frac{\widetilde{U}^{x}}{a^{x}}\right)^{n+\frac{1}{2}}\right)_{x}+\left(\widetilde{U}^{y,n+\frac{1}{2}},\psi_{x}\psi_{y}\left(\frac{\widetilde{U}^{y}}{a^{y}}\right)^{n+\frac{1}{2}}\right)_{y}\\
&+\left(\Delta t\sum_{l=0}^{n-1}\left(\frac{b^{x}}{a^{x}}\widetilde{U}^{x}\right)^{l+\frac{1}{2}},\psi_{x}\psi_{y}\left(\frac{\widetilde{U}^{x}}{a^{x}}\right)^{n+\frac{1}{2}}\right)_{x}+\left(\Delta t\sum_{l=0}^{n-1}\left(\frac{b^{y}}{a^{y}}\widetilde{U}^{y}\right)^{l+\frac{1}{2}},\psi_{x}\psi_{y}\left(\frac{\widetilde{U}^{y}}{a^{y}}\right)^{n+\frac{1}{2}}\right)_{y}\\
&+\left(\frac{\Delta t}{2}\left(\frac{b^{x}}{a^{x}}\widetilde{U}^{x}\right)^{n+\frac{1}{2}},\psi_{x}\psi_{y}\left(\frac{\widetilde{U}^{x}}{a^{x}}\right)^{n+\frac{1}{2}}\right)_{x}+\left(\frac{\Delta t}{2}\left(\frac{b^{y}}{a^{y}}\widetilde{U}^{y}\right)^{n+\frac{1}{2}},\psi_{x}\psi_{y}\left(\frac{\widetilde{U}^{y}}{a^{y}}\right)^{n+\frac{1}{2}}\right)_{y}\\
&=\left(\psi_{x}\psi_{y}f^{n+\frac{1}{2}},P^{n+\frac{1}{2}}\right).
\end{split}
\end{equation}
Similar to the derivation of Eq.$(\ref{e8})$, we have
\begin{equation}\label{e58}
  \begin{array}{r}
  \left(\partial_{t} \psi_{x}\psi_{y} P^{n+1}, P^{n+\frac{1}{2}}\right)=\frac{1}{\Delta t}\left(\left(\psi_{x}\psi_{y} P^{n+1}, P^{n+1}\right)-\left(\psi_{x}\psi_{y} P^{n}, P^{n}\right)\right)-\left(\psi_{x}\psi_{y} P^{n+\frac{1}{2}}, \partial_{t} P^{n+1}\right).
  \end{array}
  \end{equation}
Let $C_{a^{x}}=\max \left\{\left\|\frac{\partial a^{x}}{\partial x}\right\|_{\infty},\left\|\frac{\partial a^{x}}{\partial y}\right\|_{\infty}\right\},$ $C_{a^{y}}=\max \left\{\left\|\frac{\partial a^{y}}{\partial x}\right\|_{\infty},\left\|\frac{\partial a^{y}}{\partial y}\right\|_{\infty}\right\},$ and $C_{0}=min(49-6 C_{a^{x}} h_{max},49-6 C_{a^{y}} h_{max})$, we can obtain the following two estimates by using Lemma $\ref{lem3}$
\begin{equation}\label{e59}
  \begin{split}
\left(\widetilde{U}^{x,n+\frac{1}{2}},\psi_{x}\psi_{y}\left(\frac{\widetilde{U}^{x}}{a^{x}}\right)^{n+\frac{1}{2}}\right)_{x}\geq \frac{49-6 C_{a^{x}} h_{max}}{72}\left(\widetilde{U}^{x,n+\frac{1}{2}},\left(\frac{\widetilde{U}^{x}}{a^{x}}\right)^{n+\frac{1}{2}}\right)_{x}\geq \frac{C_{0}}{72a^{*}}\left\lVert \widetilde{U}^{x,n+\frac{1}{2}}\right\rVert_{x}^{2},
\end{split}
\end{equation}
\begin{equation}\label{e60}
  \begin{split}
\left(\widetilde{U}^{y,n+\frac{1}{2}},\psi_{x}\psi_{y}\left(\frac{\widetilde{U}^{y}}{a^{y}}\right)^{n+\frac{1}{2}}\right)_{y}\geq \frac{49-6 C_{a^{y}} h_{max}}{72}\left(\widetilde{U}^{y,n+\frac{1}{2}},\left(\frac{\widetilde{U}^{y}}{a^{y}}\right)^{n+\frac{1}{2}}\right)_{x}\geq \frac{C_{0}}{72a^{*}}\left\lVert \widetilde{U}^{y,n+\frac{1}{2}}\right\rVert_{y}^{2}.
\end{split}
\end{equation}
For the fourth term on the left hand side of the Eq.$(\ref{e57})$, we have that
\begin{equation}\label{e61}
  \begin{split}
    \left(\Delta t\sum_{l=0}^{n-1}\left(\frac{b^{x}}{a^{x}}\widetilde{U}^{x}\right)^{l+\frac{1}{2}},\psi_{x}\psi_{y}\left(\frac{\widetilde{U}^{x}}{a^{x}}\right)^{n+\frac{1}{2}}\right)_{x}
    &\geq -\frac{C\Delta t}{2}\sum_{l=0}^{n-1}\left\lVert \left(\frac{b^{x}}{a^{x}}\widetilde{U}^{x}\right)^{l+\frac{1}{2}}\right\rVert_{x}^{2}-\frac{\epsilon_{1}\Delta t}{2}\left\lVert \psi_{x}\psi_{y}\left(\frac{\widetilde{U}^{x}}{a^{x}}\right)^{n+\frac{1}{2}} \right\rVert_{x}^{2}\\
    &\geq -\frac{C}{2}\Delta t\sum_{l=0}^{n-1}\left\lVert \widetilde{U}^{x,l+\frac{1}{2}}\right\rVert_{x}^{2}-\frac{\epsilon_{1} \Delta t }{2\left(a_{*}\right)^{2}}\left\lVert \widetilde{U}^{x,n+\frac{1}{2}}\right\rVert_{x}^{2}.  
\end{split}
\end{equation}
Similarly,
\begin{equation}\label{e62}
  \begin{split}
    \left(\Delta t\sum_{l=0}^{n-1}\left(\frac{b^{y}}{a^{y}}\widetilde{U}^{y}\right)^{l+\frac{1}{2}},\psi_{x}\psi_{y}\left(\frac{\widetilde{U}^{y}}{a^{y}}\right)^{n+\frac{1}{2}}\right)_{y}\geq -\frac{C}{2}\Delta t\sum_{l=0}^{n-1}\left\lVert \widetilde{U}^{y,l+\frac{1}{2}}\right\rVert_{y}^{2}-\frac{\epsilon_{1} \Delta t }{2\left(a_{*}\right)^{2}} \left\lVert \widetilde{U}^{y,n+\frac{1}{2}}\right\rVert_{y}^{2}.  
\end{split}
\end{equation}
Now, we consider the sixth term on the left hand side of the Eq.$(\ref{e57})$. By using the Cauchy Schwarz inequality and Lemma $\ref{lem3},$ we have
\begin{equation}\label{e63}
  \begin{split}
    \left(\frac{\Delta t}{2}\left(\frac{b^{x}}{a^{x}}\widetilde{U}^{x}\right)^{n+\frac{1}{2}},\psi_{x}\psi_{y}\left(\frac{\widetilde{U}^{x}}{a^{x}}\right)^{n+\frac{1}{2}}\right)_{x}\leq \frac{(b^{*})^{2}+1}{2a_{*}^{2}}\Delta t \left\lVert \widetilde{U}^{x,n+\frac{1}{2}}\right\rVert_{x}^{2},
\end{split}
\end{equation}
\begin{equation}\label{e64}
  \begin{split}
    \left(\frac{\Delta t}{2}\left(\frac{b^{y}}{a^{y}}\widetilde{U}^{y}\right)^{n+\frac{1}{2}},\psi_{x}\psi_{y}\left(\frac{\widetilde{U}^{y}}{a^{y}}\right)^{n+\frac{1}{2}}\right)_{y}\leq \frac{(b^{*})^{2}+1}{2a_{*}^{2}}\Delta t \left\lVert \widetilde{U}^{y,n+\frac{1}{2}}\right\rVert_{y}^{2}.
\end{split}
\end{equation}
Then, by simple calculation, we get
\begin{equation}\label{e65}
  \begin{split}
    \left(\psi_{x}\psi_{y}f^{n+\frac{1}{2}},P^{n+\frac{1}{2}}\right)\leq C\left\lVert f^{n+\frac{1}{2}}\right\rVert^{2}+C\left\lVert P^{n+\frac{1}{2}}\right\rVert^{2}.  
\end{split}
\end{equation}
Multiplying the Eq.$(\ref{e57})$ by $2\Delta t$, summing on $n$ from 0 to $N-1$, and substituting the Eqs.$(\ref{e58})-(\ref{e65})$ into the Eq.$(\ref{e57})$, we have that
\begin{equation}\label{e66}
  \begin{split}
    &\qquad2\left(\psi_{x}\psi_{y} P^{N}, P^{N}\right)+\left(\frac{C_{0}}{36a^{*}}-\frac{(b^{*})^{2}+1}{a_{*}^{2}}\Delta t\right)\Delta t\sum_{n=0}^{N-1}\left(\left\lVert \widetilde{U}^{x,n+\frac{1}{2}}\right\rVert_{x}^{2}+\left\lVert \widetilde{U}^{y,n+\frac{1}{2}}\right\rVert_{y}^{2}\right)\\
    &\leq C\Delta t \sum_{n=0}^{N-1}\left\lVert P^{n+\frac{1}{2}}\right\rVert^{2}+\epsilon_{2}\Delta t\sum_{n=0}^{N-1}\left\lVert \partial_{t}P^{n+1}\right\rVert^{2}+C\Delta t\sum_{n=0}^{N-1}\Delta t\sum_{l=0}^{n-1}\left(\left\lVert \widetilde{U}^{x,l+\frac{1}{2}}\right\rVert_{x}^{2}+\left\lVert \widetilde{U}^{y,l+\frac{1}{2}}\right\rVert_{x}^{2}\right)\\
    &+ \frac{\epsilon_{1} \Delta t^{2} }{\left(a_{*}\right)^{2}}\sum_{n=0}^{N-1}\left(\left\lVert \widetilde{U}^{x,n+\frac{1}{2}}\right\rVert_{x}^{2}+\left\lVert \widetilde{U}^{y,n+\frac{1}{2}}\right\rVert_{y}^{2}\right)+C\Delta t\sum_{n=0}^{N-1}\left\lVert f^{n+\frac{1}{2}}\right\rVert^{2}+2\left(\psi_{x}\psi_{y} P^{0}, P^{0}\right)   
\end{split}
\end{equation} 
Taking $\epsilon_{1}$ and $\Delta t$ sufficiently small, such that
$$\frac{C_{0}}{36a^{*}}-\frac{(b^{*})^{2}+1}{a_{*}^{2}}\Delta t-\frac{\epsilon_{1} \Delta t }{a_{*}^{2}}>0,$$
and using the Gronwall's inequality and Lemma $\ref{lem2}$, we can obtain that
\begin{equation}\label{e67}
  \begin{split}
    \left\lVert P^{N}\right\rVert^{2}+\Delta t\sum_{n=0}^{N-1}\left(\left\lVert \widetilde{U}^{x,n+\frac{1}{2}}\right\rVert_{x}^{2}+\left\lVert \widetilde{U}^{y,n+\frac{1}{2}}\right\rVert_{y}^{2}\right) 
    \leq  C\Delta t\sum_{n=0}^{N}\left\lVert f^{n}\right\rVert^{2}+C\left\lVert P^{0}\right\rVert^{2} +\epsilon_{2}\Delta t\sum_{n=0}^{N-1}\left\lVert \partial_{t}P^{n+1}\right\rVert^{2}
\end{split}
\end{equation}
Next, multiplying Eq.$(\ref{e55})$ by $h_{i}k_{j}\partial_{t} P_{i,j}^{n+1}$ and summing on $i$ and $j$, we have that
\begin{equation}\label{e68}
  \begin{split}
  \left(\partial_{t}\psi_{x}\psi_{y}P^{n+1},\partial_{t} P^{n+1}\right)&+\left(U^{x,n+\frac{1}{2}},\partial_{t}\psi_{x}\psi_{y}\left(\frac{\widetilde{U}^{x}}{a^{x}}\right)^{n+1}\right)\\
  &+\left(U^{y,n+\frac{1}{2}},\partial_{t}\psi_{x}\psi_{y}\left(\frac{\widetilde{U}^{y}}{a^{y}}\right)^{n+1}\right)=\left(\psi_{x}\psi_{y}f^{n+\frac{1}{2}},\partial_{t}P^{n+1}\right).
\end{split}
\end{equation}
Then, multiplying Eq.$(\ref{e68})$ by $2\Delta t$ and summing on $n$ from $1$ to $N-1$, we obtain that
\begin{equation}\label{e69}
  \begin{split}
  2\Delta t\sum_{n=1}^{N-1}\left(\partial_{t}\psi_{x}\psi_{y}P^{n+1},\partial_{t} P^{n+1}\right)&+2\Delta t\sum_{n=1}^{N-1}\left(U^{x,n+\frac{1}{2}},\partial_{t}\psi_{x}\psi_{y}\left(\frac{\widetilde{U}^{x}}{a^{x}}\right)^{n+1}\right)\\
  &+2\Delta t\sum_{n=1}^{N-1}\left(U^{y,n+\frac{1}{2}},\partial_{t}\psi_{x}\psi_{y}\left(\frac{\widetilde{U}^{y}}{a^{y}}\right)^{n+1}\right)=2\Delta t\sum_{n=1}^{N-1}\left(\psi_{x}\psi_{y}f^{n+\frac{1}{2}},\partial_{t}P^{n+1}\right).
\end{split}
\end{equation}
By Lemma $\ref{lem2}$, it is obvious that
\begin{equation}\label{e70}
  \begin{split}
  2\Delta t\sum_{n=1}^{N-1}\left(\partial_{t}\psi_{x}\psi_{y}P^{n+1},\partial_{t} P^{n+1}\right)\geq \frac{49}{36}\sum_{n=1}^{N-1}\Delta t\left\lVert \partial_{t}P^{n+1}\right\rVert^{2} 
\end{split}
\end{equation}
Now, we consider the second term on the left hand side of the Eq.$(\ref{e69})$.
\begin{equation*}
  \begin{split}
 2\Delta t\sum_{n=1}^{N-1}\left(U^{x,n+\frac{1}{2}},\partial_{t}\psi_{x}\psi_{y}\left(\frac{\widetilde{U}^{x}}{a^{x}}\right)^{n+1}\right)&=2\Delta t\sum_{n=1}^{N-1}\left(\widetilde{U}^{x,n+\frac{1}{2}},\partial_{t}\psi_{x}\psi_{y}\left(\frac{\widetilde{U}^{x}}{a^{x}}\right)^{n+1}\right)\\
 &+2\Delta t \sum_{n=0}^{N-1}\Delta t\sum_{l=0}^{n-1}\left(\left(\frac{b^{x}}{a^{x}}\widetilde{U}^{x}\right)^{l+\frac{1}{2}},\partial_{t}\psi_{x}\psi_{y}\left(\frac{\widetilde{U}^{x}}{a^{x}}\right)^{n+1}\right)\\
 &+\Delta t\sum_{n=0}^{N-1}\Delta t\left(\left(\frac{b^{x}}{a^{x}}\widetilde{U}^{x}\right)^{n+\frac{1}{2}},\partial_{t}\psi_{x}\psi_{y}\left(\frac{\widetilde{U}^{x}}{a^{x}}\right)^{n+1}\right)\\
 &=SS_{1}+SS_{2}+SS_{3}
\end{split}
\end{equation*}
By the Lemma $\ref{lem7}$, we can easily derive the estimate of $SS_{1}$,
\begin{equation*}
  \begin{split}
   SS_{1}\geq C\left(\left\lVert \widetilde{U}^{x,N}\right\rVert_{x}^{2}-\left\lVert \widetilde{U}^{x,0}\right\rVert_{2}^{2}\right)-\epsilon_{3}\sum_{n=0}^{N-1}\Delta t\left\lVert \partial_{t}P^{n+1}\right\rVert^{2}-C\Delta t\sum_{n=0}^{N}\left\lVert \widetilde{U}^{x,n}\right\rVert_{x}^{2}.  
\end{split}
\end{equation*}
For $SS_{2}$, we have
\begin{equation*}
  \begin{split}
    SS_{2}&=2\Delta t\sum_{n=0}^{N-1}\Delta t\sum_{l=0}^{n}\left(\left(\frac{b^{x}}{a^{x}}\widetilde{U}^{x}\right)^{l+\frac{1}{2}},\psi_{x}\psi_{y}\left(\frac{\partial_{t}\widetilde{U}^{x,n+1}}{a^{x}}\right)\right)\\
    &=2\Delta t\sum_{n=1}^{N}\left(\psi_{x}\psi_{y}\left(\frac{\widetilde{U}^{x}}{a^{x}}\right)^{N}-\psi_{x}\psi_{y}\left(\frac{\widetilde{U}^{x}}{a^{x}}\right)^{n-1},\left(\frac{b^{x}}{a^{x}}\widetilde{U}^{x}\right)^{n-\frac{1}{2}}\right)\\
    &\leq 2\Delta t\sum_{n=1}^{N}\epsilon_{4}\left(\left\lVert \psi_{x}\psi_{y}\left(\frac{\widetilde{U}^{x}}{a^{x}}\right)^{N}-\psi_{x}\psi_{y}\left(\frac{\widetilde{U}^{x}}{a^{x}}\right)^{n-1}\right\rVert_{x}^{2} \right)+\frac{\Delta t}{2\epsilon_{4}}\sum_{n=1}^{N}\left\lVert \left(\frac{b^{x}}{a^{x}}\widetilde{U}^{x}\right)^{n-\frac{1}{2}}\right\rVert_{x}^{2}\\
    &\leq 2T\epsilon_{4} \left\lVert \psi_{x}\psi_{y}\left(\frac{\widetilde{U}^{x}}{a^{x}}\right)^{N}\right\rVert_{x}^{2}+2\Delta t \epsilon_{4}\sum_{n=1}^{N} \left\lVert \psi_{x}\psi_{y}\left(\frac{\widetilde{U}^{x}}{a^{x}}\right)^{n-1}\right\rVert_{x}^{2}+\frac{\Delta t}{2\epsilon_{4}}\sum_{n=1}^{N}\left\lVert \left(\frac{b^{x}}{a^{x}}\widetilde{U}^{x}\right)^{n-\frac{1}{2}}\right\rVert_{x}^{2}\\
    &\leq \frac{2T\epsilon_{4}}{a_{*}^{2}}\left\lVert \widetilde{U}^{x,N}\right\rVert_{x}^{2}+\frac{2\Delta t \epsilon_{4}}{a_{*}^{2}}\sum_{n=1}^{N}\left\lVert \widetilde{U}^{x,n-1}\right\rVert_{x}^{2}+\frac{\Delta t(b^{*})^{2}}{2\epsilon_{4} a_{*}^{2}}\sum_{n=1}^{N}\left\lVert \widetilde{U}^{x,n-\frac{1}{2}}\right\rVert_{x}^{2} \\
    &\leq  \frac{2T\epsilon_{4}}{a_{*}^{2}}\left\lVert \widetilde{U}^{x,N}\right\rVert_{x}^{2}+C\Delta t\sum_{n=0}^{N}\left\lVert \widetilde{U}^{x,n}\right\rVert_{x}^{2}. 
  \end{split}
\end{equation*}
Now, we consider $SS_{3}$,
\begin{equation*}
  \begin{split}
  SS_{3}&=\Delta t\sum_{n=0}^{N-1}\left(\left(\frac{b^{x}}{a^{x}}\widetilde{U}^{x}\right)^{n+\frac{1}{2}},\psi_{x}\psi_{y}\left(\frac{\widetilde{U}^{x,n+1}}{a^{x}}\right)-\psi_{x}\psi_{y}\left(\frac{\widetilde{U}^{x,n}}{a^{x}}\right)\right)\\
    &\leq C\Delta t\sum_{n=0}^{N-1}\left(\left\lVert \widetilde{U}^{x,n+1}\right\rVert_{x}^{2}+\left\lVert \widetilde{U}^{x,n}\right\rVert_{x}^{2}\right)\\
    &\leq C\Delta t\sum_{n=0}^{N}\left\lVert \widetilde{U}^{x,n}\right\rVert_{x}^{2}.
  \end{split}
\end{equation*}
Therefore,
\begin{equation}\label{e71}
  \begin{split}
 2\Delta t\sum_{n=1}^{N-1}\left(U^{x,n+\frac{1}{2}},\partial_{t}\psi_{x}\psi_{y}\left(\frac{\widetilde{U}^{x}}{a^{x}}\right)^{n+1}\right)\geq &\left(C-\frac{2T\epsilon_{4}}{a_{*}^{2}}\right)\left\lVert \widetilde{U}^{x,N}\right\rVert_{x}^{2}-C\left\lVert \widetilde{U}^{x,0}\right\rVert_{x}^{2}\\
 &-\epsilon_{3}\sum_{n=0}^{N-1}\Delta t\left\lVert \partial_{t}P^{n+1}\right\rVert^{2}-C\Delta t\sum_{n=0}^{N}\left\lVert \widetilde{U}^{x,n}\right\rVert_{x}^{2}.
\end{split}
\end{equation}
Similarly,
\begin{equation}\label{e72}
  \begin{split}
 2\Delta t\sum_{n=1}^{N-1}\left(U^{y,n+\frac{1}{2}},\partial_{t}\psi_{x}\psi_{y}\left(\frac{\widetilde{U}^{y}}{a^{y}}\right)^{n+1}\right)\geq &\left(C-\frac{2T\epsilon_{4}}{a_{*}^{2}}\right)\left\lVert \widetilde{U}^{y,N}\right\rVert_{y}^{2}-C\left\lVert \widetilde{U}^{y,0}\right\rVert_{y}^{2}\\
 &-\epsilon_{3}\sum_{n=0}^{N-1}\Delta t\left\lVert \partial_{t}P^{n+1}\right\rVert^{2}-C\Delta t\sum_{n=0}^{N}\left\lVert \widetilde{U}^{y,n}\right\rVert_{y}^{2}.
\end{split}
\end{equation}
Next, we estimate the term on the right hand side of the Eq.$(\ref{e69})$,
\begin{equation}\label{e73}
  \begin{split}
    2\Delta t\sum_{n=1}^{N-1}\left(\psi_{x}\psi_{y}f^{n+\frac{1}{2}},\partial_{t}P^{n+1}\right)\leq C\Delta t\sum_{n=0}^{N}\left\lVert \psi_{x}\psi_{y}f^{n}\right\rVert^{2}+\epsilon_{5}\Delta t\sum_{n=0}^{N-1}\left\lVert \partial_{t}P^{n+1}\right\rVert^{2}.
\end{split}
\end{equation}
Taking $\epsilon_{3}, \epsilon_{4}$ and $\epsilon_{5}$ sufficiently small such that 
$$\frac{49}{36}-2\epsilon_{3}-\epsilon_{5}>0,$$
$$C-\frac{2T\epsilon_{4}}{a_{*}^{2}}>0,$$
and combining the Eqs.$(\ref{e69})-(\ref{e73})$, we have that
\begin{equation}\label{e74}
  \begin{split}
    \left\lVert \widetilde{U}^{x,N}\right\rVert_{x}^{2}+\left\lVert \widetilde{U}^{y,N}\right\rVert_{y}^{2}+\left(\frac{49}{36}-2\epsilon_{3}-\epsilon_{5}\right)\Delta t\sum_{n=0}^{N-1}\left\lVert \partial_{t}P^{n+1}\right\rVert^{2}\leq &C\left\lVert \widetilde{U}^{x,0}\right\rVert_{x}^{2}
    +C\Delta t\sum_{n=0}^{N}\left\lVert \widetilde{U}^{x,n}\right\rVert_{x}^{2}+C\left\lVert \widetilde{U}^{y,0}\right\rVert_{y}^{2}\\
    &+C\Delta t\sum_{n=0}^{N}\left\lVert \widetilde{U}^{y,n}\right\rVert_{y}^{2}+C\Delta t\sum_{n=0}^{N}\left\lVert \psi_{x}\psi_{y}f^{n}\right\rVert^{2}.
\end{split}
\end{equation}
Using the Gronwall's inequality, we obtain that
\begin{equation}\label{e75}
  \begin{split}
    \left\lVert \widetilde{U}^{x,N}\right\rVert_{x}^{2}+\left\lVert \widetilde{U}^{y,N}\right\rVert_{y}^{2}&+\left(\frac{49}{36}-2\epsilon_{3}-\epsilon_{5}\right)\Delta t\sum_{n=0}^{N-1}\left\lVert \partial_{t}P^{n+1}\right\rVert^{2}\\
    &\leq C\left\lVert \widetilde{U}^{x,0}\right\rVert_{x}^{2}
    +C\left\lVert \widetilde{U}^{y,0}\right\rVert_{y}^{2}+C\Delta t\sum_{n=0}^{N}\left\lVert f^{n}\right\rVert^{2}.
\end{split}
\end{equation}
Taking $\epsilon_{2}$ sufficiently small such that
$$\frac{49}{36}-2\epsilon_{3}-\epsilon_{5}-\epsilon_{2}>0,$$
and combining the Eq.$(\ref{e67})$ and Eq.$(\ref{e75}),$ we have
\begin{equation}\label{e76}
  \begin{split}
    \left\lVert P^{N}\right\rVert^{2}+\left\lVert \widetilde{U}^{x,N}\right\rVert_{x}^{2}+\left\lVert \widetilde{U}^{y,N}\right\rVert_{y}^{2}\leq C\left(\left\lVert p^{0}\right\rVert^{2} +\left\lVert \widetilde{U}^{x,0}\right\rVert_{x}^{2}+\left\lVert \widetilde{U}^{y,0}\right\rVert_{y}^{2}+\Delta t\sum_{n=0}^{N}\left\lVert f^{n}\right\rVert^{2}\right).
\end{split}
\end{equation}
\begin{thm}\label{thm3}
  Suppose $h$ and $\Delta t$ are sufficiently small. Let $P$ and $\widetilde{U}$ be the solutions to CBCFD scheme. There exists a positive constant $C$ such that the following inequality holds,
  \begin{equation}\label{e77}
    \begin{split}
      \left\lVert P^{N}\right\rVert^{2}+\left\lVert \widetilde{U}^{x,N}\right\rVert_{x}^{2}+\left\lVert \widetilde{U}^{y,N}\right\rVert_{y}^{2}\leq C\left(\left\lVert p^{0}\right\rVert^{2} +\left\lVert \widetilde{U}^{x,0}\right\rVert_{x}^{2}+\left\lVert \widetilde{U}^{y,0}\right\rVert_{y}^{2}+\Delta t\sum_{n=0}^{N}\left\lVert f^{n}\right\rVert^{2}\right).
  \end{split}
  \end{equation}
\end{thm}
\subsection{The error estimate}
In this section, we will discuss the convergence of the CBCFD scheme. Suppose $p, u^{x},u^{y}\in C^{3}(0,T;C^{5}(\Omega))$. From the Eq.$(\ref{e49})$, we can get the following equations by Taylor's expansion.
\begin{equation}\label{e78}
  \begin{split}
  \partial_{t} \psi_{x}\psi_{y} p_{i,j}^{n+1}+\delta_{x} u_{i,j}^{x,n+\frac{1}{2}}+\delta_{y} u_{i,j}^{y,n+\frac{1}{2}}=\psi_{x}\psi_{y} f_{i,j}^{n+\frac{1}{2}}+O\left(h^{4}+k^{4}+\Delta t^{2}\right),
  \end{split}
  \end{equation}
  \begin{equation}\label{e79}
    \begin{split}
    \delta_{x} p_{i+1/2,j}^{n+1}=-\psi_{x}\left(\frac{\widetilde{u}^{x}}{a^{x}}\right)_{i+1/2,j}^{n+1}+O\left(h^{4}\right),
    \end{split}
  \end{equation}
  \begin{equation}\label{e80}
    \begin{split}
    \delta_{y} p_{i,j+1/2}^{n+1}=-\psi_{y}\left(\frac{\widetilde{u}^{y}}{a^{y}}\right)_{i,j+1/2}^{n+1}+O\left(k^{4}\right),
    \end{split}
  \end{equation}
Set
\begin{equation*}
    \begin{split}
    \eta&=P-p,\\
    \widetilde{\xi}&=\widetilde{U}^{x}-\widetilde{u}^{x},\\
    \xi&=U^{x}-u^{x}=U^{x}-\widetilde{U}^{x}+\widetilde{\xi}^{x}+\widetilde{u}^{x}-u^{x},\\ 
    \widetilde{\gamma}&=\widetilde{U}^{y}-\widetilde{u}^{y},\\
    \gamma&=U^{y}-u^{y}=U^{y}-\widetilde{U}^{y}+\widetilde{\xi}^{y}+\widetilde{u}^{y}-u^{y}.
    \end{split}
\end{equation*}
Subtracting Eq.$(\ref{e78})$ from Eq.$(\ref{e55})$, Eq.$(\ref{e79})$ from Eq.$(\ref{e51})$ and Eq.$(\ref{e80})$ from Eq.$(\ref{e52})$, we have
\begin{equation}\label{e81}
  \begin{split}
    \partial_{t} \psi_{x}\psi_{y} \eta_{i,j}^{n+1}+\delta_{x} \xi_{i,j}^{n+\frac{1}{2}}+\delta_{y} \gamma_{i,j}^{n+\frac{1}{2}}=R_{1,i,j}^{n+1},
  \end{split}
  \end{equation}
  \begin{equation}\label{e82}
    \begin{split}
    \delta_{x} \eta_{i+1/2,j}^{n+1}=-\psi_{x}\left(\frac{\widetilde{\xi}}{a^{x}}\right)_{i+1/2,j}^{n+1}+R_{2,i+\frac{1}{2},j}^{n+1},
    \end{split}
  \end{equation}
  \begin{equation}\label{e83}
    \begin{split}
    \delta_{y} \eta_{i,j+1/2}^{n+1}=-\psi_{y}\left(\frac{\widetilde{\gamma}}{a^{y}}\right)_{i,j+1/2}^{n+1}+R_{3,i,j+1/2}^{n+1},
    \end{split}
  \end{equation}
  where $R_{1,i}^{n+1}=O(h^{4}+k^{4}+\Delta t^{2}),$ $R_{2,i+\frac{1}{2},j}^{n+1}=O(h^{4}),$ and $R_{3,i,j+1/2}^{n+1}=O(k^{4})$.

  The convergence result can be obtained by a similar process to the Theorems $\ref{thm2}$ and $\ref{thm3}$. 
  \begin{thm}\label{thm4}
    Suppose $h$ and $\Delta t$ are sufficiently small and $p, u^{x},u^{y}\in C^{3}(0,T;C^{5}(\Omega))$. Let $P$, $\widetilde{U}^{x}$ and $\widetilde{U}^{y}$ be the solutions to CBCFD scheme. There exists a positive constant $C$ such that
    \begin{equation}\label{e22}
      \begin{split}
        \left\lVert P^{N}-p^{N}\right\rVert+\left\lVert \widetilde{U}^{x,N}-\widetilde{u}^{x,N}\right\rVert+\left\lVert \widetilde{U}^{y,N}-\widetilde{u}^{y,N}\right\rVert\leq C\left(\Delta t^2+h^{4}+k^{4} \right).  
    \end{split}
    \end{equation}
  \end{thm}
\section{Numerical experiments}
To verify the validity of the CBCFD scheme, we will carry out two examples including the exact solutions of polynomial functions and trigonometric types. The Example 1 is an one-dimensional case with the domain $\Omega=(0,1)$. The Example 2 is a two-dimensional model with the domain $\Omega=(0,1)\times (0,1)$. In the two examples, we set $J=(0,1]$ and take $\Delta t=h^{2}$ for showing the fourth spatial convergence order. The uniform grids are available in this section. In addition, we compute the results of the normal BCFD scheme and display the comparison of errors for the two methods which indicates that the CBCFD method has a higher accuracy than the normal CBCFD method. 

\begin{figure}
  \centering
  \subfloat{%
  \resizebox*{10cm}{!}{\includegraphics{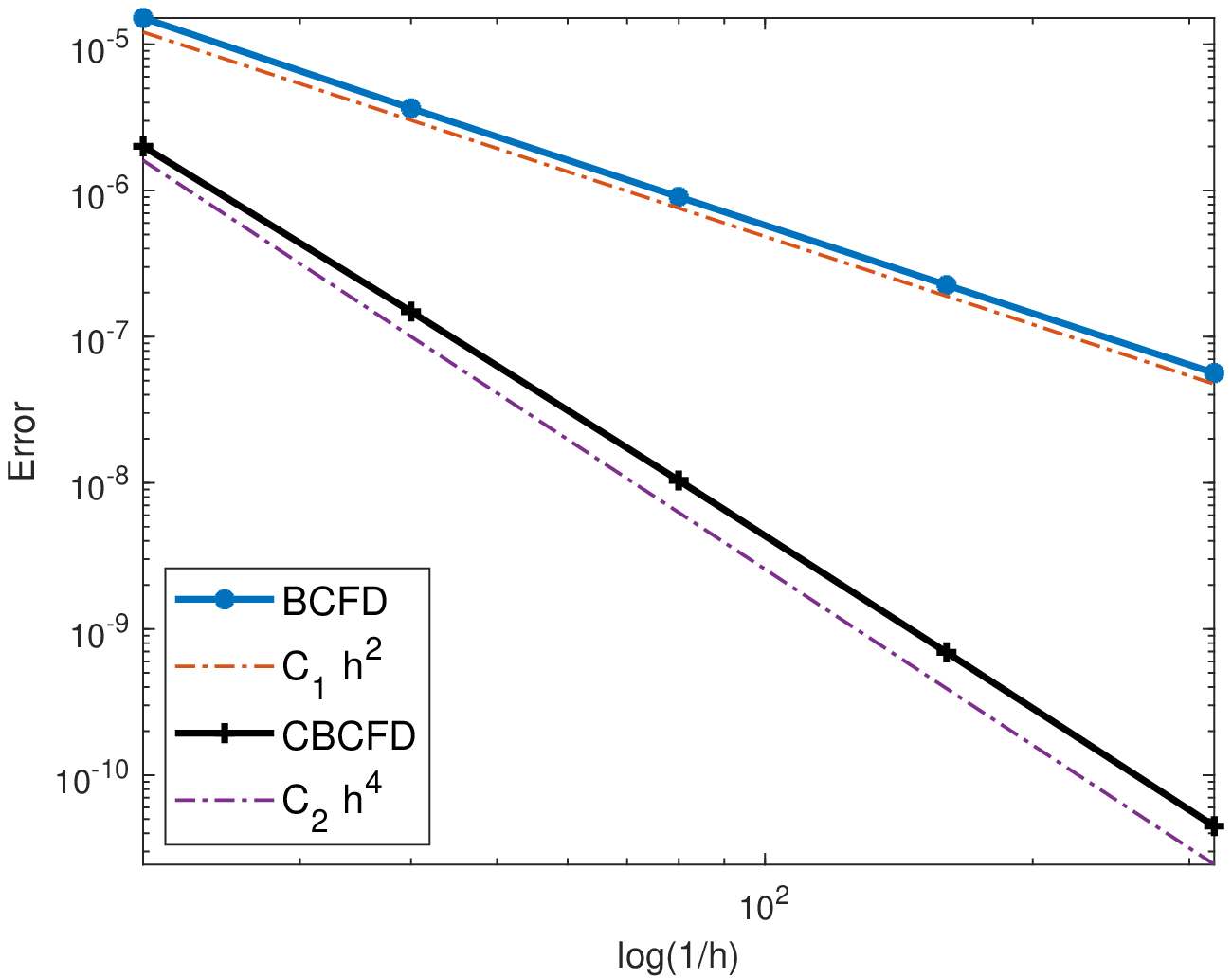}}}\hspace{5pt}
  \caption{The comparison of convergence orders of the two schemes \\
  \hspace{-1.5cm}BCFD and CBCFD in the Example 1.} \label{fig:example1p2}
\end{figure}

\textbf{Example 1:} In this experiment, the numerical results for an one-dimensional case will be displayed. We take the coefficients $a$ and $b$ and the exact solution $p$ as follows. The source/sink term $f$ can be obtained by the direct calculation. The errors, convergence orders and the comparison between BCFD and CBCFD scheme are listed in Table \ref{tab:example1} and Figure \ref{fig:example1p2}.
\begin{equation*}
  \left\{\begin{array}{ll}  p=tx^4(1-x)^4, \\
    a=1.0\times 10^{-8},\\
    b=1,\\
    f=x^4(1-x)^4-4.0\times 10^{-8}x^{2}(1-x)^{2}(3-14x+11x^{2}-\frac{3}{2}t^{2}+7xt^{2}+3x^{2}t-7x^{2}t^{2}).
\end{array}\right.
\end{equation*}

 \begin{table}[htbp]
  \centering
  \setlength{\tabcolsep}{6.5mm}
  \caption{Error and convergence rates in $h$ of example 1.}
  \label{tab:example1}
  \begin{tabular}{lcccccc} 
  \hline
   $h$ & $\|p^{N}-Z^{N}\|_{L^{2}(\Omega)}$ & $Rates$ & $\|\widetilde{u}^{N}-\widetilde{U}^{N}\|_{L^{2}(\Omega)}$ & $Rates$ \\ 
  \hline
  1/20  & 2.00E-06 & ---   & 2.89E-14 & ---   \\
  1/40  & 1.48E-07 & 3.7553  & 2.19E-15 & 3.7218   \\
  1/80  & 1.04E-08 & 3.8377  & 1.52E-16 & 3.8541   \\
  1/160  & 6.92E-10 & 3.9072  & 9.96E-18 & 3.9286    \\
  1/320  & 4.47E-11 & 3.9506  & 6.38E-19 & 3.9654   \\
  \hline
  \end{tabular}
\end{table}

\begin{figure}
  \centering
  \subfloat{%
  \resizebox*{10cm}{!}{\includegraphics{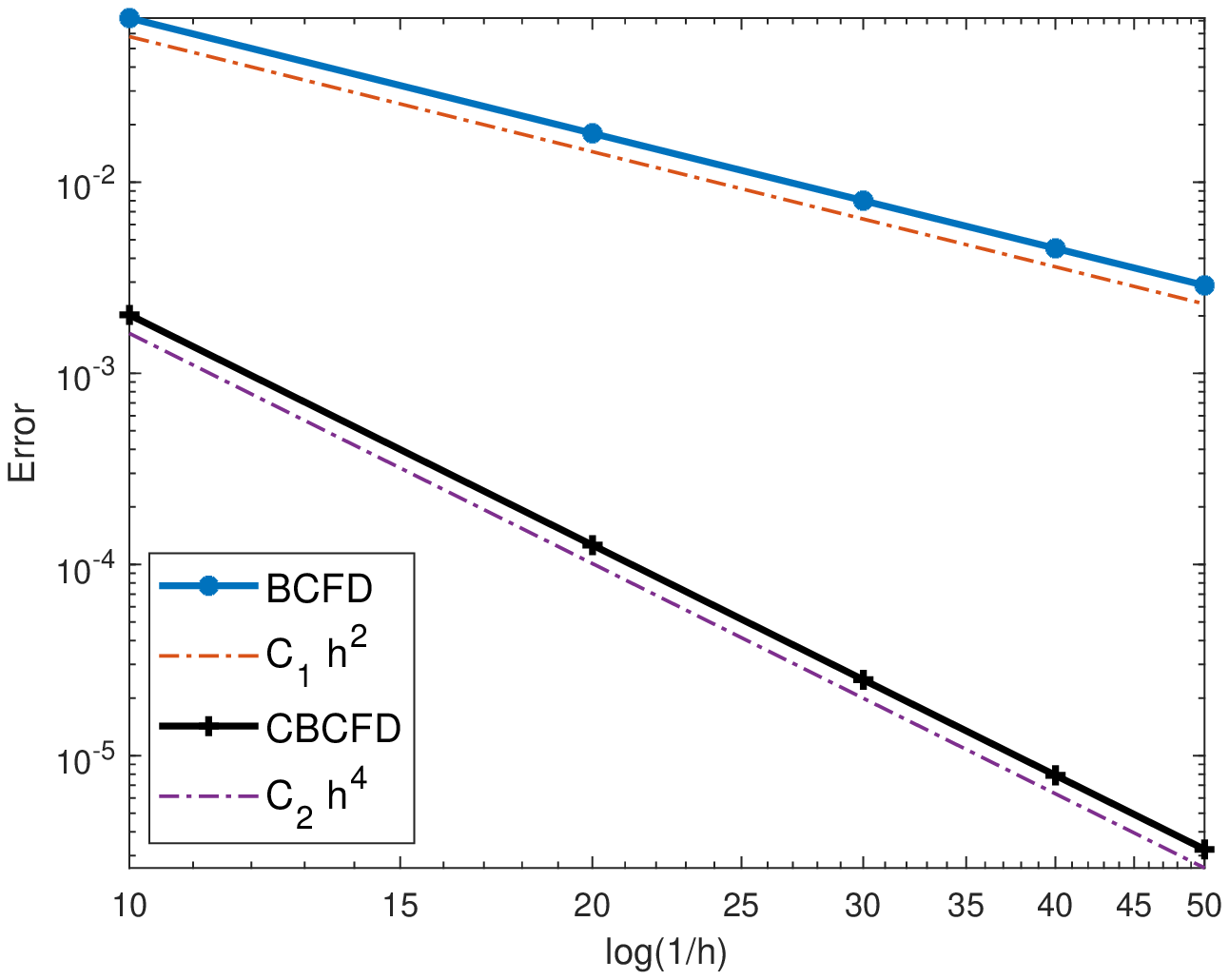}}}\hspace{5pt}
  \caption{The comparison of convergence orders of the two schemes \\
  \hspace{-1.5cm}BCFD and CBCFD in the Example 2.} \label{fig:example2p2}
\end{figure}

 \textbf{Example 2:} We consider a two-dimensional case and take the coefficient matrices $\textbf{A}$ and $\textbf{B},$ the exact solution $p$ and the source/sink term $f$ as follows. The numerical results are listed in Table \ref{tab:example2} and Figure \ref{fig:example2p2}.
 \begin{equation*}
  \left\{\begin{array}{ll}  p=t^{2}*cos(2\pi x)*cos(2\pi y), \\
    \textbf{A}=\textbf{I}, \\
    \textbf{B}=t\textbf{I},\\
    f=(2t+8\pi^{2}t^{2}+\frac{\pi^{2}}{2}t^{4})cos(2\pi x)cos(2\pi y).
\end{array}\right.
\end{equation*}

\begin{table}[htbp]
  \centering
  \setlength{\tabcolsep}{6.5mm}
  \caption{Error and convergence rates in $h$ of example 2.}
  \label{tab:example2}
  \begin{tabular}{lcccccc} 
  \hline
   $h$ & $\|p^{N}-Z^{N}\|_{L^{2}(\Omega)}$ & $Rates$ & $\|\widetilde{u}^{N}-\widetilde{U}^{N}\|_{L^{2}(\Omega)}$ & $Rates$ \\ 
  \hline
  1/10  & 4.52E-04 & ---   & 1.97E-03 & ---   \\
  1/20  & 2.82E-05 & 4.0016  & 1.23E-04 & 4.0011   \\
  1/30  & 5.58E-06 & 4.0004  & 2.43E-05 & 4.0003   \\
  1/40  & 1.76E-06 & 4.0002  & 7.69E-06 & 4.0002    \\
  1/50  & 7.23E-07 & 4.0001  & 3.15E-06 & 4.0001   \\
  \hline
  \end{tabular}
\end{table}
 
From the data in Table \ref{tab:example1}-\ref{tab:example2} and Figure \ref{fig:example1p2}-\ref{fig:example2p2}, we observe that the numerical solutions approximate the exact results well. In addition, we can conclude that the CBCFD schemes are valid for solving the non-Fickian flow models and have the  fourth convergence order in spatial direction.
 \section{Conclusion} 
In this work, the numerical schemes combined the compact block-centered finite difference methods in spatial direction with the Crank-Nicolson discretization in temporal  direction were constructed and analyzed 
to solve non-Fickian flow models in porous media. We also established the stability analyses and error estimates of the constructed schemes rigorously both in one-dimensional and two-dimensional cases. Finally we carried out some numerical experiments to verify the theoretical analysis.

 \section*{Acknowledgments}
 This work was supported by National Natural Science Foundation of China (12131014,11901489).



\bibliographystyle{elsarticle-num}
\bibliography{manuscript}

\newpage

\end{document}